\documentclass{aims}
\usepackage{paralist}
\usepackage{bm}
\usepackage{amsmath, amssymb, graphicx}

\usepackage{mathtools}

\usepackage{cite}
\usepackage{setspace}
 \usepackage[colorlinks=true]{hyperref}

\def\currentvolume{X}
 \def\currentissue{X}
  \def\currentyear{200X}
   \def\currentmonth{XX}
    \def\ppages{X--XX}
     \def\DOI{10.3934/xx.xx.xx.xx}
\renewcommand{\eqref}[1]{(\ref{#1})}
\newcommand{\mat}[1]{\bm{#1}}
\renewcommand{\vec}[1]{\bm{#1}}
\newcommand{\manifold}{\ensuremath{\mathcal{M}}}
\newcommand{\map}{\ensuremath{\vec F}}
\newcommand{\koopman}{\ensuremath{\mathcal{K}}}
\newcommand{\observables}{\ensuremath{\mathcal{F}}}
\newcommand{\dictionary}{\ensuremath{\mathcal{D}}}
\newcommand{\dmd}{\ensuremath{\mat{K}_{\text{DMD}}}}
\newcommand{\expect}{\mathbb{E}}
\newcommand{\dmdvec}{\vec{\xi}}
\newcommand{\dmdmat}{\mat{\Xi}}

\begin{document}

\title[Approximating the Koopman Operator with Data] 
      {A Data--Driven Approximation \\of the Koopman Operator:\\
        Extending  Dynamic Mode Decomposition}
\author[M.O. Williams \and I.G. Kevrekidis \and C.W. Rowley]{}

 \email{mow2@princeton.edu}
 \email{krevrekidis@princeton.edu}
 \email{cwrowley@princeton.edu}


\begin{abstract}

The Koopman operator is a {\em linear} but infinite dimensional operator that governs the evolution of scalar observables defined on the state space of an autonomous dynamical system, and is a powerful tool for the analysis and decomposition of nonlinear dynamical systems.
In this manuscript, we present a data driven method for approximating the leading {\em eigenvalues, eigenfunctions, and modes} of the Koopman operator. 
The method requires a data set of snapshot pairs and a dictionary of scalar observables, but does not require explicit governing equations or interaction with a ``black box'' integrator.
We will show that this approach is, in effect, an extension of Dynamic Mode Decomposition (DMD), which has been used to approximate the Koopman eigenvalues and modes.
Furthermore, if the data provided to the method are generated by a Markov process instead of a deterministic dynamical system, the algorithm approximates the eigenfunctions of the Kolmogorov backward equation, which could be considered as the ``stochastic Koopman operator''~\cite{Mezic2005}.
Finally, four illustrative examples are presented: two that highlight the quantitative performance of the method when presented with either deterministic or stochastic data, and two that show potential applications of the Koopman eigenfunctions. 
\end{abstract}

\maketitle

\centerline{\scshape Matthew O. Williams }
\medskip
{\footnotesize
 \centerline{Program in Applied and Computational Mathematics (PACM)}
   \centerline{Princeton University, NJ 08544, USA.}
} 

\medskip
\centerline{\scshape Ioannis G. Kevrekidis }
\medskip
{\footnotesize
 \centerline{Chemical and Biological Engineering Department \& PACM}
   \centerline{Princeton University, NJ 08544, USA.}
} 

\medskip
\centerline{\scshape Clarence W. Rowley }
\medskip
{\footnotesize
 \centerline{Mechanical and Aerospace Engineering Department}
   \centerline{Princeton University, NJ 08544, USA.}
} 

\bigskip


\section{Introduction}
\label{sec:intro}

In many mathematical and engineering applications, a phenomenon of interest can be summarized in different ways. 
For instance, to describe the state of a two dimensional incompressible fluid flow, one can either record velocity and pressure fields or streamfunction and vorticity~\cite{Hirsch2007}.
Furthermore, these states can often be {\em approximated} using a low dimensional set of Proper Orthogonal Decomposition (POD) modes~\cite{Holmes1998}, a set of Dynamic Modes~\cite{Schmid2010,Schmid2012}, or a finite collection of Lagrangian particles~\cite{Monaghan1992}.
A mathematical example is the linear time invariant (LTI) system provided by $\vec x(n+1) = \mat{A} \vec x(n)$, where $\vec x(n)$ is the system state at the $n$-th timestep.
Written as such, the evolution of $\vec x$ is governed by the eigenvalues of $\mat{A}$. 
One could also consider the invertible but nonlinear transformation, $\vec z(n) = \vec T(\vec x(n))$, which generates a nonlinear evolution law for   $\vec z$.
Both  approaches (i.e., $\vec x$ or $\vec z$) describe the same fundamental behavior, yet one description may be preferable to others.
For example, solving an LTI system is almost certainly preferable to evolving a nonlinear system from a computational standpoint.

In general, one measures (or computes) the state of a system using a set of scalar {\em observables}, which are functions defined on state space, and watches how the values of these functions evolve in time.
Furthermore, provided the set of observations is rich enough, one can even write an evolution law for the dynamics of the set of observations, and use this system in lieu of the original one.
Because the properties of this new dynamical system depend on our choice of variables (observables), it would be highly desirable if one could find a set of observables whose dynamics appear to be governed by a linear evolution law.
If such a set could be identified, the dynamics would be completely determined by the spectrum of the evolution operator.
Furthermore, this could enable the simple yet effective algorithms designed for linear systems, for example controller design~\cite{Todorov2006,Stengel2012} or stability analysis~\cite{Mezic2005,Lehoucq1998}, to be applied to nonlinear systems. 

Mathematically, the evolution of observables of the system state is governed by the {\em Koopman operator}~\cite{Koopman1932,Koopman1931,Budivsic2012,Rowley2009}, which is a {\em linear but infinite dimensional} operator that is  defined for an autonomous dynamical system.
Of particular interest here is the ``slow'' subspace of the Koopman operator, which is the span of the eigenfunctions associated with eigenvalues near the unit circle in discrete time (or near the imaginary axis in continuous time).
These eigenvalues and eigenfunctions capture the long term dynamics of observables that appear after the fast transients have subsided, and could serve as a low dimensional approximation of the otherwise infinite dimensional operator when a spectral gap, which clearly delineates the ``fast'' and ``slow'' temporal dynamics, is present.
In addition to the eigenvalues and eigenfunctions, the final element of Koopman spectral analysis is the set of {\em Koopman modes} for the  {\em full state} observable~\cite{Rowley2009, Budivsic2012}, which are vectors that enable us to reconstruct the \textit{state of the system} as a linear combination of the Koopman eigenfunctions.
Overall, the ``tuples'' of Koopman eigenfunctions, eigenvalues, and modes enable us to: (a) transform state space so that it the dynamics appear to be linear, (b) determine the temporal dynamics of the linear system, and (c) reconstruct the state of the original system from our new linear representation.
In principle, this framework is quite broadly applicable, and useful even for problems with multiple attractors that cannot be accurately approximated using models based on local linearization.

There are several algorithms in the literature that can computationally approximate subsets of these quantities. 
Three examples are Generalized Laplace Analysis (GLA)~\cite{Budivsic2012,Mauroy2012,Mauroy2013}, the Ulam Galerkin Method~\cite{Froyland2013,Bollt2013}, and Dynamic Mode Decomposition (DMD)~\cite{Schmid2010, Tu2013, Rowley2009}.
None of these techniques require explicit governing equations, so all, in principle, can be applied directly to data. 
GLA can approximate both the Koopman modes and eigenfunctions, but it requires knowledge of the eigenvalues to do so~\cite{Budivsic2012,Mauroy2012,Mauroy2013}.
The Ulam Galerkin method has been used to approximate the eigenfunctions and eigenvalues~\cite{Froyland2013}, though it is more frequently used to generate finite dimensional approximations of the Perron--Frobenius operator, which is the adjoint of the Koopman operator.
Finally, DMD has been used to approximate the Koopman modes and eigenvalues~\cite{Rowley2009,Tu2013}, but not the Koopman eigenfunctions.

Even in pairs instead of triplets, approximations of these quantities are useful. 
DMD and its variants~\cite{Wynn2013,Chen2012,Jovanovic2013} have been successfully used to analyze nonlinear fluid flows using data from both experiments and computation~\cite{Schmid2010,Muld2012,Seena2011}. 
GLA and similar methods have been applied to extract meaningful spatio-temporal structures using sensor data from buildings and power systems \cite{Eisenhower2010,Susuki2011,Susuki2012,Susuki2013}.
Finally, the Ulam Galerkin method has been used to identify coherent structures and almost invariant sets~\cite{Froyland2007,Froyland2009,Froyland2005} based on the singular value decomposition of  (a slight  modification of) the  Koopman operator.

In this manuscript, we present a data driven method that approximates the leading Koopman {\em eigenfunctions, eigenvalues, {\bf and} modes} from a data set  of successive ``snapshot'' pairs and a dictionary of observables that spans a subspace of the scalar observables.
There are many possible ways to choose this dictionary, and it could be comprised of polynomials, Fourier modes, spectral elements, or other sets of functions of the full state observable.
We will argue that this approach is an extension of DMD that can produce better approximations of the Koopman eigenfunctions; as such, we refer to it as \emph{Extended Dynamic Mode Decomposition} (EDMD). 
One regime where the behavior of both EDMD and DMD can be formally analyzed and contrasted is in the limit of large data.
In this regime, we will show that the numerical approximation of the Koopman eigenfunctions generated by EDMD converges to the numerical approximation we would obtain from a Galerkin method~\cite{Boyd2013} in that the residual is orthogonal to the subspace spanned by the elements of the dictionary.
With finite amounts of data, we will demonstrate the effectiveness of EDMD on two deterministic examples: one that highlights the quantitative accuracy of the method, and a more practical  application.

Because EDMD is an entirely data driven procedure, it can also be applied to data from stochastic systems without any algorithmic changes.
If the underlying system is a Markov process, we will show that EDMD approximates the eigenfunctions of the Kolmogorov backward equation~\cite{Givon2004, Bagherisubmitted}, which has been called the  \textit{stochastic  Koopman} operator (SKO)~\cite{Mezic2005}. 
Once again, we will demonstrate the effectiveness of the EDMD procedure when the amount of data is limited by applying it to two stochastic examples: the first to test the accuracy of the method, and the second to highlight a potential application of EDMD as a nonlinear manifold learning technique. 
In the latter example, we highlight two forms of model reduction: reduction that occurs when the \textit{dynamics of the system state} are constrained to a low-dimensional manifold, and reduction that occurs when the {\em statistical moments} of the stochastic dynamical system are effectively low dimensional.

In the remainder of the manuscript, we will detail the EDMD algorithm and show (when mathematically possible) or demonstrate through examples that it accurately approximates the leading Koopman eigenfunctions, eigenvalues, and modes for both deterministic and stochastic sets of data.
In particular, in Sec.~\ref{sec:dmd} the EDMD algorithm will be presented, and we will prove that it converges to a Galerkin approximation of the Koopman operator given a sufficiently large amount of data. 
In Sec.~\ref{sec:dictionary}, we detail three choices of dictionary that we have found to be effective in a broad set of applications.
In Sec.~\ref{sec:deterministic}, we will demonstrate that the EDMD approximation can be accurate even with finite amounts of data, and can yield useful parameterizations of common dynamical structures such as systems with multiple basins of attraction {\em when the underlying system is deterministic}.
In Sec.~\ref{sec:stochastic-ck}, we experiment by applying EDMD to stochastic data and show it approximates the eigenfunctions of the SKO for Markov processes. 
Though the interpretation of the eigenfunctions now differs, we demonstrate that they can still be used to accomplish useful tasks such as the parameterization of nonlinear manifolds. 
Finally, some brief concluding remarks are given in Sec.~\ref{sec:conclusions}.

\section{Dynamic Mode Decomposition and the Koopman Operator}
\label{sec:dmd}

Our ambition in this section is to establish the connection between the Koopman operator and what we call EDMD.
To accomplish this, we will define the Koopman operator in Sec.~\ref{subsec:koopman}.
Using this definition, we will outline the EDMD algorithm in Sec.~\ref{subsec:edmd}, and then show how it can be used to approximate the Koopman eigenvalues, eigenfunctions, and modes.
Next, in Sec.~\ref{subsec:convergence}, we will prove that the EDMD method almost surely converges to a Galerkin method in the limit of large data.
Finally, in Sec.~\ref{subsec:dmd}, we will highlight the connection between the EDMD algorithm and standard DMD.

\subsection{The Koopman Operator}
\label{subsec:koopman}

\begin{figure}
\centering
\includegraphics[width=\textwidth]{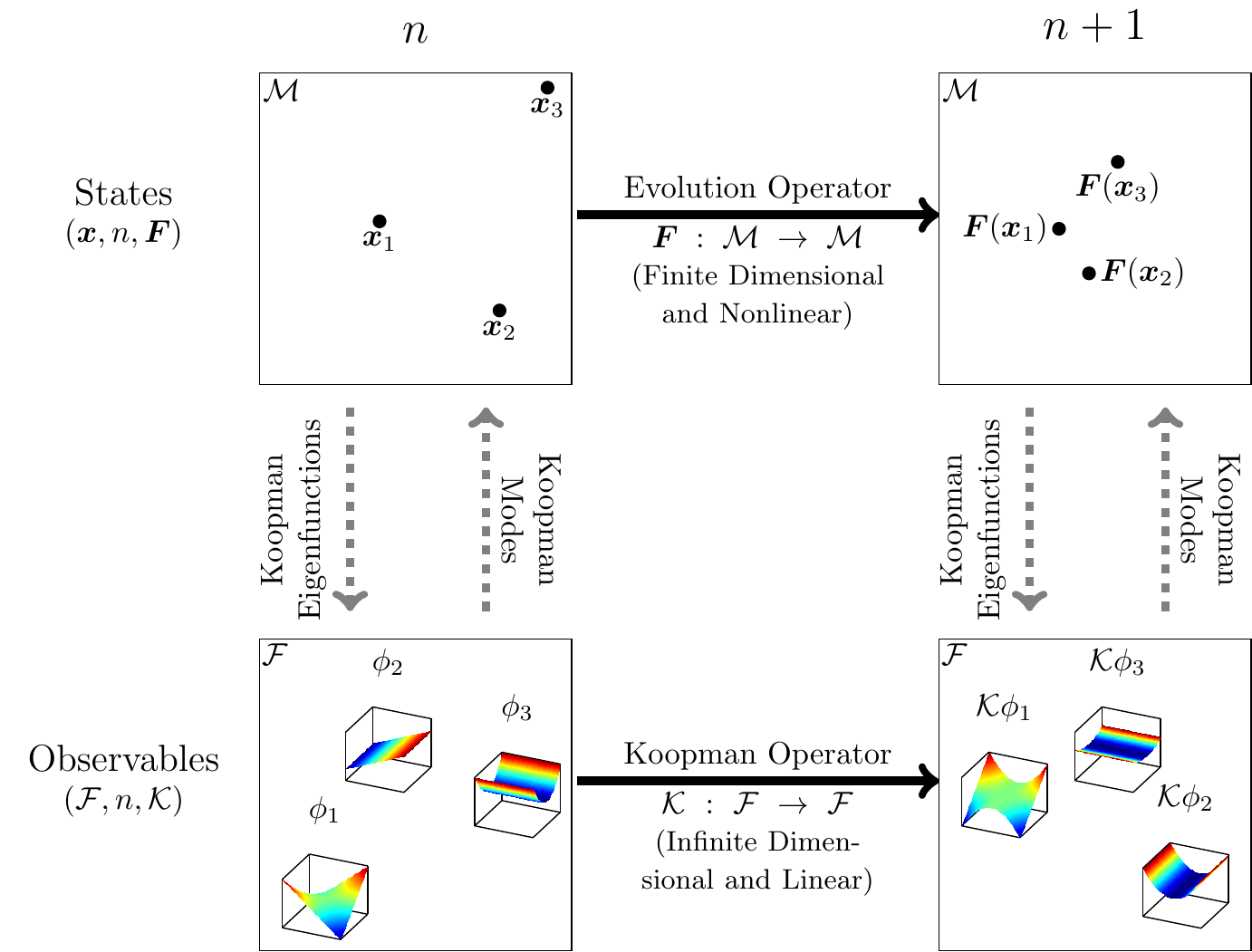}
\caption{ 
A cartoon of the Koopman operator and how it relates to the underlying dynamical system.
The top path updates the state, $\vec x\in\manifold$, using the evolution operator $\map$.
The bottom path updates the observables, $\phi\in\observables$, using the Koopman operator, $\koopman$.
Here both dynamical systems are autonomous, so the (discrete) time, $n\in\mathbb{Z}$, does not appear explicitly.
The connection between the states and observables is through the {\em full state observable} $\vec g(\vec x) = \vec x$. 
By writing $\vec g$ in terms of the Koopman eigenfunctions, we substitute the complex evolution of $\vec x$ with the straightforward, linear evolution of the $\{\varphi_i\}_i$.
To reconstruct $\vec x$, we superimpose the Koopman eigenfunctions evaluated at a point, which satisfy $(\koopman\varphi_i)(\vec{ x}_i) = \mu_i\varphi_i(\vec x_i)$, using the Koopman modes as shown in \eqref{eq:koopman-predict}.
As a result, these two ``paths'' commute, and one can either solve a finite dimensional but nonlinear problem (the top path) or an infinite dimensional but linear problem (the bottom path) if one can compute the Koopman eigenvalues, eigenfunctions, and modes.
}
\label{fig:comparison}
\end{figure}

Because the Koopman operator is central to all that follows, we will define it along with the properties relevant to EDMD in this subsection. 
No new mathematical results are presented here; our only objective  is to include for completeness the terms and definitions we will require later in the paper.
To do so, we require the autonomous, discrete time dynamical system $(\manifold, n, \map)$, where $\manifold\subseteq \mathbb{R}^N$ is the state space, $n\in\mathbb{Z}$ is (discrete) time, and $\map:\manifold\to\manifold$ is the evolution operator.
Unlike $\map$, which acts on $\vec x \in\manifold$, the Koopman operator, $\koopman$, acts on functions of state space, $\phi\in\observables$ with $\phi:\manifold\to\mathbb{C}$.
The action of the Koopman operator is,
\begin{equation}
\koopman \phi = \phi \circ \map, 
\label{eq:koopman}
\end{equation}
where $\circ$ denotes the composition of $\phi$ with $\map$.   
We stress once again that {\em the Koopman operator maps functions of state space to functions of state space} and not states to states~\cite{Koopman1931,Koopman1932,Budivsic2012}.

In essence, the Koopman operator defines a new dynamical system, $(\observables, n, \koopman)$, that governs the evolution of {\em observables}, $\phi\in\observables$, in discrete time. 
In what follows, we assume that $\observables = L^2(\manifold, \rho)$, where $\rho$ is a positive, single valued analytic function defined on $\manifold$, but not necessarily an invariant measure of the underlying dynamical system.
This assumption, which has been made before in the literature~\cite{Budivsic2012,Koopman1931}, is required so that the inner products in the Galerkin-like method we will present can be taken.
Because it acts on functions, $\koopman$ is infinite dimensional even when $\map$ is finite dimensional, but it is also {\em linear even when $\map$ is nonlinear}. 
The infinite dimensional nature of the Koopman operator is potentially problematic,
but if it can, practically, be truncated without too great a loss of accuracy (e.g., if the system has multiple time scales), then the result would be a linear and finite dimensional approximation. 
Therefore, the promise of the Koopman approach is to take the tools developed for linear systems and apply them to the dynamical system defined by the Koopman operator; thus obtaining a linear approximation of a nonlinear system without directly linearizing around a fixed point.

The dynamical system defined by $\map$ and the one defined by $\koopman$ are two different parameterizations of the same fundamental behavior.
The link between these parameterizations is the ``full state observable,'' $\vec g(\vec x) = \vec x$ and $\{(\mu_k, \varphi_k, \vec v_k)\}_{k=1}^K$, the set of $K$ ``tuples'' of Koopman eigenvalues, eigenfunctions, and modes required to reconstruct the full state. 
Note that  $K$ could (and often will) be infinite.
Although $\vec g$ is a {\em vector valued observable}, each component of it is a scalar valued observable, i.e., $g_i\in\observables$ where $g_i$ is the $i$-th component of $\vec g$. 
Assuming $g_i$ is in the span of our set of $K$ eigenfunctions, $g_i = \sum_{k=1}^K v_{ik}\varphi_k$ with $v_{kj}\in\mathbb{C}$. 
Then $\vec g$ can be obtained by ``stacking'' these weights into vectors (i.e., $\vec v_j = [v_{1j}, v_{2j},\ldots, v_{Nj}]^T$). 
As a result, 
\begin{equation}
\vec x = \vec g(\vec x) = \sum_{k=1}^K \vec v_k \varphi_k(\vec x),
\label{eq:koopman-modes}
\end{equation}
where $\vec v_k$ is the $k$-th {\em Koopman mode}, and $\varphi_k$ is the $k$-th {\em Koopman eigenfunction}.
In doing this, we have assumed that each of the scalar observables that comprise $\vec g$ are in the subspace of $\observables$ spanned by our $K$ eigenfunctions, but we have not assumed that the eigenfunctions form a basis for $\observables$.

The system state at future times can be obtained either by directly evolving $\vec x$ or by evolving the full state observable through Koopman:  
\begin{align}
\map(\vec x) &= \left(\koopman \vec g\right)(\vec x)  
= \sum_{k=1}^K \vec v_k (\koopman\varphi_k)(\vec x) = \sum_{k=1}^K \mu_k \vec v_k \varphi_k(\vec x).
\label{eq:koopman-predict}
\end{align}
This representation of $\map(\vec x)$ is particularly advantageous because the dynamics associated with each eigenfunction are determined by its corresponding eigenvalue. 

Figure~\ref{fig:comparison} shows a commutative diagram that acts as a visual summary of this section. 
The top row shows the direct evolution of states, $\vec x \in \manifold$, governed by $\map$; the bottom row shows the evolution of observables, $\phi \in \observables$, governed by the Koopman operator. 
Although $\map$ and $\koopman$ act on different spaces, they encapsulate the same dynamics. 
For example, once given a state $\vec x$, to compute $(\koopman\phi)(\vec x)$ one could either take the observable $\phi$, apply $\koopman$, and evaluate it at $\vec x$ (the bottom route), or use $\map$ to compute $\map(\vec x)$ and then evaluate $\phi$ at this updated position (the top route).
Similarly, to compute $\map(\vec x)$, one could either apply $\map$ to $\vec x$ (the top route) or apply $\koopman$ to the full state observable and evaluate $(\koopman\vec g)(\vec x)$ (the bottom route).
As a result, one can either choose to work with a finite dimensional, nonlinear system or an infinite dimensional, linear system depending upon which ``path'' is simpler/more useful for a given problem.

\subsection{Extended Dynamic Mode Decomposition}
\label{subsec:edmd}

In this subsection, we outline Extended Dynamic Mode Decomposition (EDMD), which is a method that approximates the Koopman operator and therefore the Koopman eigenvalue, eigenfunction, and mode tuples defined in Sec.~\ref{subsec:koopman}.
The EDMD procedure requires: (a) a data set of snapshot pairs, i.e., $\{(\vec x_m, \vec y_m)\}_{m=1}^M$ that we will organize as {\em a pair of data sets},
\begin{equation}
\mat{X} = 
\begin{bmatrix}
\vec x_1 & \vec x_2 & \cdots & \vec x_M
\end{bmatrix}, 
\qquad 
\mat{Y} = 
\begin{bmatrix}
\vec y_1 & \vec y_2 & \cdots & \vec y_M
\end{bmatrix},
\label{eq:data}
\end{equation}
where $\vec x_i\in\manifold$ and $\vec y_i\in\manifold$ are {\em snapshots of the system state} with  $\vec y_i = \map(\vec x_i)$, and  (b) a dictionary of observables, $\dictionary=\{\psi_1, \psi_2, \ldots, \psi_K\}$ where $\psi_i\in\observables$, whose span we denote as $\observables_{\dictionary}\subset\observables$; for brevity, we also define the vector valued function $\vec{\Psi}:\manifold\to\mathbb{C}^{1\times K}$ where 
\begin{equation}
\vec\Psi(\vec x) = 
\begin{bmatrix}
\psi_1(\vec x) & \psi_2(\vec x) & \cdots &  \psi_K(\vec x)
\end{bmatrix}
. 
\label{eq:function-def}
\end{equation}
The data set needed is typically constructed from multiple short bursts of simulation or from experimental data. 
For example, if the data were given as a single time series, then for a given snapshot $\vec x_i$, $\vec y_i = \map(\vec x_i)$ is the next snapshot in the time series. 
The optimal choice of dictionary elements  remains an open question, but a short discussion including some pragmatic choices will be given in Sec.~\ref{sec:dictionary}.
For now, we assume that $\dictionary$ is ``rich enough'' to accurately approximate a few of the leading Koopman eigenfunctions.

\subsubsection{Approximating the Koopman Operator and its Eigenfunctions}
Now we seek to generate $\mat{K}\in\mathbb{R}^{K\times K}$, a finite dimensional approximation of $\koopman$.
By definition, a function $\phi\in\observables_{\dictionary}$ can be written as 
\begin{equation}
\phi = \sum_{k=1}^K a_k\psi_k = \vec\Psi\vec a,
\label{eq:function-approx}
\end{equation}
the linear superposition of the $K$ elements in the dictionary with the weights $\vec a$.
Because $\observables_{\dictionary}$ is typically not an invariant subspace of $\koopman$,
\begin{equation}
\koopman\phi = (\vec\Psi\circ\map) \vec a = \vec\Psi(\mat{K}\vec a) + r
\label{eq:koopman-residual}
\end{equation}
which includes the residual term $r\in\observables$.
To determine $\mat{K}$, we will minimize
\begin{equation}
\label{eq:objective-function}
\begin{aligned}
J &= \frac{1}{2}\sum_{m=1}^M |r(\vec x_m)|^2 \\
&= \frac{1}{2}\sum_{m=1}^M \left|((\vec\Psi\circ \map)(\vec x_m) - \vec\Psi(\vec x_m)\mat{K})\vec a \right|^2 \\
&= \frac{1}{2}\sum_{m=1}^M \left|(\vec\Psi(\vec y_m) - \vec\Psi(\vec x_m)\mat{K})\vec a \right|^2 \\
\end{aligned}
\end{equation} 
where $\vec x_m$ is the $m$-th snapshot in $\mat{X}$, and $\vec y_m = \vec F(\vec x_m)$ is the $m$-th snapshot in $\mat{Y}$. 
Equation~\ref{eq:objective-function} is a least squares problem, and therefore cannot have multiple isolated local minima; it must either have a unique global minimizer or a continuous family (or families) of minimizers. 
As a result, regularization (here via the truncated singular value decomposition) may be required to ensure the solution is unique, and the $\mat{K}$ that minimizes (\ref{eq:objective-function}) is:
\begin{equation}
\mat{K} \triangleq \mat{G}^+\mat{A}, 
\label{eq:edmd}
\end{equation}
where $+$ denotes the pseudoinverse and 
\begin{align}
\mat{G} &= \frac{1}{M}\sum_{m=1}^M \vec\Psi(\vec x_m)^* \vec\Psi(\vec x_m),
\\ 
\mat{A} &= \frac{1}{M}\sum_{m=1}^M \vec\Psi(\vec x_m)^* \vec\Psi(\vec y_m),
\end{align}
with $\mat{K},\mat{G},\mat{A}\in\mathbb{C}^{K\times K}$. 
As a result, $\mat{K}$ is a finite dimensional approximation of $\koopman$ that maps $\phi\in\observables_{\dictionary}$ to some other  $\hat\phi\in\observables_{\dictionary}$ by minimizing the residuals at the data points.
As a consequence, if $\dmdvec_j$ is the $j$-th eigenvector of $\mat{K}$ with the eigenvalue $\mu_j$, then the EDMD approximation of an eigenfunction of $\koopman$ is
\begin{equation}
\varphi_j = \vec\Psi\dmdvec_j.
\label{eq:edmd-eigenfunction}
\end{equation}
Finally, in many applications  the discrete time data in $\mat{X}$ and $\mat{Y}$ are generated by a continuous time process with a sampling interval of $\Delta t$.
If this is the case, we define $\lambda_j = \frac{\ln(\mu_j)}{\Delta t}$ to approximate the eigenvalues of the \textit{continuous time system}.
In the remainder of the manuscript, we denote the eigenvalues of $\mat{K}$ with the $\mu_j$ and (when applicable) the approximation of the corresponding continuous time eigenvalues as $\lambda_j$. 
Although both embody the same information, one choice is often more natural for a specific problem.

\subsubsection{Computing the Koopman Modes}
\label{subsec:koopman-modes}

Next, we will compute approximations of the Koopman modes for the full state observable using EDMD.
Recall that the Koopman modes are the weights needed to express the full state in the {\em Koopman eigenfunction} basis.
As such, we will proceed in two steps: first, we will express the full state observable using the elements of $\dictionary$; then, we will find a mapping from the elements of $\dictionary$ to the numerically computed eigenfunctions.
Applying these two steps in sequence will yield the observables expressed as a linear combination of Koopman eigenfunctions, which is, by definition, the Koopman modes for the full state observable.

Recall that the full state observable, $\vec g(\vec x) = \vec x$, is a vector valued observable (e.g., $\vec g:\manifold\to\mathbb{R}^N$) that can be generated by ``stacking'' $N$ scalar valued observables, $g_i:\manifold\to\mathbb{R}$, as follows: 
\begin{equation}
\vec g(\vec x) = 
\begin{bmatrix}
g_1(\vec x)\\
g_2(\vec x) \\
\vdots \\
g_N(\vec x)
\end{bmatrix}
=
\begin{bmatrix}
\vec e_1^*\vec x \\
\vec e_2^*\vec x \\
\vdots \\
\vec e_N^* \vec x
\end{bmatrix},
\end{equation}
where $\vec e_i$ is the $i$-th unit vector in $\mathbb{R}^N$.
At this time, we conveniently assume that all $g_i\in\observables_{\dictionary}$ so that $g_i = \sum_{k=1}^K \psi_k b_{k,i} = \vec\Psi\vec b_i$, where $\vec b_i$ is some appropriate vector of weights.
If this is not the case, {\em approximate} Koopman modes can be computed by projecting $g_i$ onto $\observables_\dictionary$, though the accuracy and usefulness of this fit clearly depends on the choice of $\dictionary$.
To avoid this issue,  $g_i\in\observables_\dictionary$ for $i=1,\ldots,N$ in all  examples that follow. 
In either case, the entire vector valued observable can be expressed (or approximated) in this manner as 
\begin{equation}
\vec g = \mat{B}^T 
\mat{\Psi}^T = (\vec\Psi\mat{B})^T, \qquad 
\mat{B} = 
\begin{bmatrix}
\vec b_1 & \vec b_2 & \cdots & \vec b_N
\end{bmatrix},
\label{eq:vector-valued-basis}
\end{equation}
where $\mat{B}\in\mathbb{C}^{K\times N}$.

Next, we will express the $\psi_i$ in terms of all the $\varphi_i$, which are our numerical approximations of the Koopman eigenfunctions.
For notational convenience, we define the vector--valued functions $\vec\Phi:\manifold\to\mathbb{C}^K$, where
\begin{equation}
\label{eq:eigenfunction-vector}
\mat{\Phi}(\vec x) = 
\begin{bmatrix}
\varphi_1(\vec x) & \varphi_2(\vec x) & \cdots & \varphi_K(\vec x)
\end{bmatrix}.
\end{equation}
Using \eqref{eq:edmd-eigenfunction} and \eqref{eq:vector-valued-basis}, this function can also be written as 
\begin{equation}
\mat{\Phi}
= \mat{\Psi}\mat{\Xi}
,
\qquad
\dmdmat = 
\begin{bmatrix}
\dmdvec_1 & \dmdvec_2 & \cdots & \dmdvec_K
\end{bmatrix},
\end{equation}
where $\dmdvec_i\in\mathbb{C}^K$ is the $i$-th eigenvector of $\mat{K}$ associated with $\mu_i$. 
Therefore, we can determine the $\psi_i$ as a function of $\varphi_i$ by inverting $\dmdmat^T$. 
Because $\dmdmat$ is a matrix of eigenvectors, its inverse is 
\begin{equation}
\dmdmat^{-1} = \mat{W}^* = 
\begin{bmatrix}
\vec w_1 & \vec w_2 & \cdots \vec w_K
\end{bmatrix}^*,
\label{eq:inverse}
\end{equation}
where $\vec w_i$ is the $i$-th left eigenvector of $\mat{K}$ also associated with $\mu_i$ (i.e., $\vec w_i^*\mat{K} = \vec w_i^*\mu_i$) appropriately scaled so $\vec w_i^*\dmdvec_i = 1$.
We combine \eqref{eq:vector-valued-basis} and \eqref{eq:inverse}, and after some slight algebraic manipulation find that
\begin{equation}
\vec g = 
\mat{V}
\mat{\Phi}^T
= 
\sum_{k=1}^K \vec v_k \varphi_k, \qquad
\mat{V} = 
\begin{bmatrix}
\vec v_1 & \vec v_2 & \ldots & \vec v_K
\end{bmatrix}
=\left(\mat{W}^*\mat{B}\right)^T,
\label{eq:compute-modes}
\end{equation}
where $\vec v_i = (\vec w_i^*\mat{B})^T$ is the $i$-th Koopman mode.
This is the formula for the Koopman modes that we desired.

In summary, EDMD requires a data set of snapshot pairs, $\{(\vec x_m, \vec y_m)\}_{m=1}^M$, that we represent as two data sets, $\mat{X}$ and $\mat{Y}$, as well as a dictionary of observables, $\dictionary$.
Furthermore, it assumes that the leading Koopman eigenfunctions are (nearly) contained within $\observables_\dictionary$, the subspace spanned by the elements of  $\dictionary$. 
With this information, a finite dimensional approximation of the Koopman operator, $\mat{K}$, can be computed using \eqref{eq:edmd}. 
The eigenvalues of $\mat{K}$ are the EDMD approximations of the Koopman eigenvalues.
The {\em right eigenvectors} of $\mat{K}$ generate the approximations of the eigenfunctions, and the {\em left eigenvectors} of $\mat{K}$ generate the approximations of the Koopman modes.

\subsection{Convergence of the EDMD Algorithm to a Galerkin Method}
\label{subsec:convergence}

In this subsection, we relate EDMD to the Galerkin methods one would use to approximate the Koopman operator with complete information about the underlying dynamical system.
In this context, a Galerkin method is a weighted residual method where the residual, as defined in \eqref{eq:koopman-residual}, is orthogonal to the span of $\dictionary$. 
In particular, we show that the EDMD approximation of the Koopman operator converges to the approximation that would be obtained from a Galerkin method when $M$ becomes sufficiently large and if: (1) the elements of $\mat{X}$ are drawn from a distribution on $\manifold$ with density $\rho$, and (2) $\observables=L^2(\manifold, \rho)$.
The first assumption  defines a process for adding new data points to our set, and could be replaced with other sampling schemes.
The second assumption is required so that the inner products in the Galerkin method converge, which is relevant  for problems where  $\manifold=\mathbb{R}^N$.

If EDMD were a Galerkin method, then the entries of $\mat{G}$ and $\mat{A}$ in (\ref{eq:edmd}) would be defined as 
\begin{equation}
\label{eq:dmd-galerkin}
\begin{aligned}
 \hat{\mat G}_{ij} &= \int_{\manifold} \psi_i^*(\vec x)\psi_j(\vec x)\rho(\vec x)\;d\vec x = 
        \left\langle\psi_i, \psi_j\right\rangle_\rho, \\ 
 \hat{\mat A}_{ij} &= \int_{\manifold} \psi_i^*(\vec x)\psi_j(\map(\vec x))\rho(\vec x)\;d\vec x = 
\left\langle\psi_i, \koopman\psi_j\right\rangle_\rho,
\end{aligned}
\end{equation}
where $\left\langle p,q\right\rangle_\rho = \int_{\manifold} p^*(\vec x)q(\vec x)\rho(\vec x) \; d\vec x$ is the inner product, and the finite dimensional Galerkin approximation of the Koopman operator would be $\hat{\mat{K}} = \hat{\mat{G}}^{-1}\hat{\mat{A}}$. 
The performance of this  method certainly depends upon the choice of $\psi_j$ and $\rho$, but it is nevertheless a Galerkin method as the residual would be orthogonal to $\observables_\dictionary$~\cite{Boyd2013,Trefethen2000}. 
There are nontrivial questions about what sets of $\psi$ and what measures, $\rho$, are required if the Galerkin method is to generate a {\em useful} approximation of the Koopman operator (e.g., when can we ``trust'' our eigenfunctions if $\rho$ is compactly supported but $\manifold=\mathbb{R}^N$?), but they are beyond the scope of this manuscript and will be the focus of future work.

For a finite $M$, the $ij$-th element of $\mat{G}$ is 
\begin{subequations}
\label{eq:dmd-finite}
\begin{equation}
\mat G_{ij} \triangleq \frac{1}{M}\sum_{m=1}^M\psi_i^*(\vec x_m)\psi_j(\vec x_m) = \overline{\psi_i ^*(\vec x)\psi_j(\vec x)},
\end{equation}
where the bar denotes the sample mean.
Similarly, 
\begin{equation}
\mat A_{ij} \triangleq \frac{1}{M}\sum_{m=1}^M\psi_i^*(\vec x_m)\psi_j(\vec y_m) = \overline{\psi_i^*(\vec x) (\psi_j\circ\map)(\vec x)}.
\end{equation}
\end{subequations}
When $M$ is finite, \eqref{eq:dmd-galerkin} is approximated by \eqref{eq:dmd-finite}. 
However by the law of large numbers, the sample means almost surely converge to the expectation when the number of samples, $M$, becomes sufficiently large. 
For this system, the expectation can be written as 
\begin{equation}
\begin{aligned}
\lim_{M\to\infty} \mat{G}_{ij} &= \int_{\manifold} \psi_i^*(\vec x)\psi_j(\vec x)\rho(\vec x)\;d\vec x = 
        \left\langle\psi_i, \psi_j\right\rangle_\rho = \hat{\mat{G}}_{ij}, \\ 
\lim_{M\to\infty} \mat{A}_{ij} &= \int_{\manifold} \psi_i^*(\vec x)\psi_j(\map(\vec x))\rho(\vec x)\;d\vec x = 
\left\langle\psi_i, \koopman\psi_j\right\rangle_\rho = \hat{\mat{A}}_{ij},
\end{aligned}
\end{equation}
which reintroduces the integrals in \eqref{eq:dmd-galerkin}.
As a result, the entries of $\mat{A}$ and $\mat{G}$ converge to the values they would have if the integrals were taken analytically, and therefore, the output of the EDMD procedure will converge to the output of a Galerkin method.
With randomly distributed initial data, the needed integrals are computed using Monte-Carlo integration, and the rate of convergence will be $\mathcal{O}(M^{-1/2})$.
Other sampling choices, such as placing points on a uniform grid, effectively use different quadrature rules and could therefore obtain a better rate of convergence.

\subsection{Relationship with DMD}
\label{subsec:dmd}

When $M$ is not large, EDMD will not be an accurate Galerkin method because the quadrature errors generated by the Monte-Carlo integrator will be significant, and so the residual will probably not  be orthogonal to $\observables_{\dictionary}$. 
However, it is still formally an extension of DMD, which has empirically been shown to yield meaningful results even without exhaustive data sets. 
In this section, we show that EDMD is equivalent to DMD {\em for a very specific -- and restrictive -- choice of $\dictionary$} because  EDMD and DMD will produce the same set of eigenvalues and modes for any set of snapshot pairs.
Because there are many conceptually equivalent but mathematically different definitions of DMD, the one we adopt here is taken from Ref.~\cite{Tu2013}, which defines the {\em DMD modes} as the eigenvectors of the matrix
\begin{equation}
\dmd \triangleq \mat{Y}\mat{X}^+,
\label{eq:exact-dmd}
\end{equation}
where the $j$-th mode is associated with the $j$-th eigenvalue of $\dmd$, $\mu_j$.
$\dmd$ is constructed using the data matrices in \eqref{eq:data},  where $+$ again denotes the pseudoinverse.
This definition is a generalization of preexisting DMD algorithms~\cite{Schmid2010,Schmid2011}, and does not require the data to be in the form of a single time series.

Now, we will prove that the Koopman modes computed using EDMD with $\dictionary=\{\vec e_1^*, \vec e_2^*, \ldots, \vec e_N^*\}$, {\em which is the special (if relatively restrictive) choice of dictionary alluded to earlier}, are equivalent to the DMD modes by showing that the $i$-th Koopman mode, $\vec v_i$, is also an eigenvector of $\dmd$ and, hence, a DMD mode.
Because the elements of the full state observable are the dictionary elements, $\mat{B}=\mat{I}$ in \eqref{eq:vector-valued-basis}.
Then,  the Koopman modes are the complex conjugates of the left eigenvectors of $\mat{K}$, so $\vec v_i^T = \vec w_i^*$.
Furthermore,  $\mat{G}^T = \frac{1}{M}\mat{X}\mat{X}^*$ and $\mat{A}^T = \frac{1}{M}\mat{Y}\mat{X}^*$.
Then 
\begin{equation}
\mat{K}^T = \mat{A}^T\mat{G}^{T+} = \mat{Y}\mat{X}^*\left(\mat{X}\mat{X}^*\right)^+
= \mat{Y}\mat{X}^+ = \dmd.
\end{equation}
Therefore, $\dmd\vec v_i = (\vec v_i^T\dmd^T)^T = (\vec w_i^*\mat{K})^T = (\mu_i \vec w_i^*)^T = \mu_i \vec v_i$, and all the Koopman modes computed by EDMD are eigenvectors of $\dmd$ and, thus, the DMD modes.  
Once again, the choice of dictionary is critical; {\em EDMD and DMD are only equivalent for this very specific $\dictionary$}, and other choices of $\dictionary$ will generate different (and potentially more useful) results.

Conceptually,  DMD can be thought of as {\em producing an approximation of the Koopman eigenfunctions} using the set of linear monomials as basis functions for $\observables_{\dictionary}$,  which  is analogous to a one--term Taylor expansion.
For problems where the eigenfunctions can be approximated accurately using linear monomials (e.g., in some small neighborhood of a stable fixed point), then DMD will produce an {\em accurate}  local approximation of the Koopman eigenfunctions. 
However, this is certainly not the case for all systems (particularly beyond the region of validity for local linearization). 
EDMD can be thought of as an extension of DMD that retains additional terms in the expansion, where these additional terms are determined by the elements of $\dictionary$. 
The quality of the resulting approximation is governed by $\observables_{\dictionary}$, and therefore, depends upon the choice of $\dictionary$.
However, the hope is that a more extensive $\dictionary$ will produce a superior approximation of the Koopman eigenfunctions compared to the one produced by DMD simply it uses a larger $\observables_{\dictionary}$.
As a result, with the right choice of $\dictionary$, EDMD should be applicable to a broader array of problems where the implicit choice of $\dictionary$ made by DMD results in less accuracy than desired even if $M$ is not large enough for EDMD to be an accurate Galerkin method.

\section{The Choice of the Dictionary}
\label{sec:dictionary}

As with all spectral methods, the accuracy and rate of convergence of EDMD depends on $\dictionary$, whose elements, which we refer to as trial or basis functions, span the subspace of observables, $\observables_\dictionary\subset \observables$. 
Possible choices for the elements of $\dictionary$ include: polynomials~\cite{Boyd2013}, Fourier modes~\cite{Trefethen2000}, radial basis functions~\cite{Wendland1999}, and spectral elements ~\cite{Karniadakis2013}, but the optimal choice of basis functions likely depends on both the underlying dynamical system and the sampling strategy used to obtain the data. 
Any of these sets are, in principle, a useful choice for $\dictionary$, though some care must be taken on infinite domains to ensure that any needed inner products will converge.

Choosing $\dictionary$ for EDMD is, in some cases, more difficult than selecting a set of basis functions for use in a standard spectral method because the domain on which the underlying dynamical system is defined, $\manifold$, is not necessarily known.
Typically, we can define $\Omega\supset\manifold$ so that it contains all the data in $\mat{X}$ and $\mat{Y}$, e.g., pick $\Omega$ to be a ``box'' in $\mathbb{R}^N$ that contains every snapshot in $\mat{X}$ and $\mat{Y}$.  
Next, we choose the elements of $\dictionary$ to be a basis for $\tilde\observables_{\dictionary}\subset\tilde \observables$ where $\tilde\observables$ is the space of functions that map $\Omega\to\mathbb{C}$.
Because $\observables\subset\tilde\observables$, this choice of $\dictionary$ can be used in the EDMD procedure, but there is no guarantee that the elements of $\dictionary$ form a basis for $\observables_{\dictionary}$ as there may be redundancies. 
The potential for these redundancies and the numerical issues they generate is why regularization and hence the pseudoinverse~\cite{Hansen1990} is required in \eqref{eq:edmd}.
An example of these redundancies and their effects is given in App.~\ref{app:redundant}.

\begin{table}
\centering
\caption{A table with some commonly used sets of trial functions, and the application where they are, on our experience, most suited.
}

\begin{tabular}{lp{2.5in}}
\textbf{Name} &  \textbf{Suggested Context}  \\ 
\hline \\
Hermite Polynomials & Problems defined on $\mathbb{R}^N$  \\ 
Radial Basis Functions & Problems defined on irregular domains  \\ 
Discontinuous Spectral Elements & Large problems where a block-diagonal $\mat{G}$ is beneficial/computationally important     \\ 
\end{tabular}
\label{table:basis-choices} 
\end{table} 
 
Although the optimal choice of $\dictionary$ is unknown, there are three choices that are broadly applicable in our experience.
They are: Hermite polynomials, radial basis functions (RBFs), and discontinuous spectral elements.
The \textit{Hermite polynomials} are the simplest of the three sets, and are best suited to problems defined on $\mathbb{R}^N$ if the data in $\mat{X}$ are normally distributed. 
The observables that comprise $\dictionary$ are products of the Hermite polynomials in a single dimension  (e.g., $H_1(x)H_2(y)H_0(z)$ where $H_i$ is the $i$-th Hermite polynomial and $\vec x = (x,y,z)$). 
This set of basis functions is simple to implement, and conceptually related to approximating the Koopman eigenfunctions with a Taylor expansion.
Furthermore, because they are orthogonal with respect to Gaussian weights, $\mat{G}$ will be diagonal if the $\vec x_m$ are drawn from a normal distribution, which can be beneficial numerically.

An alternative to the Hermite polynomials are \textit{discontinuous spectral elements}. 
To use this set, we define a set of $B_N$ boxes, $\{\mathcal{B}_i\}_{i=1}^{B_N}$, such that $\cup_{i=1}^{B_N}\mathcal{B}_i \supset\manifold$. 
Then, on each of the $\mathcal{B}_i$, we define $K_i$ (suitably transformed) Legendre polynomials. 
For example, in one dimension, each basis function is of the form
\begin{equation}
\psi_{ij}(x) = 
\begin{cases}
L_j(\xi) & x\in\mathcal{B}_i, \\
0          & \text{ otherwise},
\end{cases}
\end{equation}
where $L_j$ is the $j$-th Legendre polynomial, and $\xi$ is $x$ transformed such that the ``edges'' of the box are at $\xi = \pm 1$.
The advantage of this basis is that $\mat G$ will be block diagonal, and therefore easy to invert even if a very large number of basis functions are employed.

With a fixed amount of data, an equally difficult task is choosing the $\mathcal{B}_i$; the number and arrangement of the $\mathcal{B}_i$ is a balance between span of the basis functions (i.e., $h$-type convergence), which increases as the number of boxes is increased, and the accuracy of the quadrature rule, which decreases because smaller boxes contain fewer data points. 
To generate a covering of $\manifold$, we use a  method similar to the one used by GAIO~\cite{Dellnitz2001}. 
Initially, all the data (i.e., $\mat{X}$ and $\mat{Y}$) are contained within a single user selected box, $\mathcal{B}_1^{(0)}$.
If this box contains more than a pre-specified number of data points, it is subdivided into $2^N$ domains of equation measure (e.g., in one dimension, $\mathcal{B}_1^{(0)} = \mathcal{B}_1^{(1)} \cup \mathcal{B}_2^{(1)}$). 
We then proceed recursively, if any of $\mathcal{B}_i^{(1)}$ contain more than a pre-specified number of points, then they too are subdivided; this proceeds until no box has an ``overly large'' number of data points.
Any $\mathcal{B}_i^{(j)}$ that do not contain any data points are pruned, which after $j$ iterates leaves the set of subdomains, $\{\mathcal{B}_i^{(j)}\}$, on which we define the Legendre polynomials.  
The resulting trial functions are compactly supported and can be evaluated efficiently using $2^N$ trees, where $N$ is the dimension of a snapshot.
Finally, the higher order polynomials used here allow for more rapid $p$-type convergence if the eigenfunctions happen to be smooth.

The final choice of trial functions is a set of \textit{radial basis functions} (RBFs), which appeal to previous work on  ``mesh-free'' methods~\cite{Liu2010}.
Because these methods do not require a computational grid or mesh, they are particularly effective for problems where $\manifold$ has what might be called a complex geometry.
Many different RBFs could be effective, but one particularly useful set of RBFs are the thin plate splines~\cite{Wendland1999,Belytschko1996} because they do not require the scaling parameter that other RBFs (e.g., Gaussians) do.  
However, we still must choose the ``centers'' about which the RBFs are defined, which we do with $k$-means clustering~\cite{Bishop2006} with a pre-specified value of $k$ on the combined data set.
Although we make no claims of optimality, in our examples, the density of the RBF centers appears to be directly related to the density of data points, which is, intuitively, a reasonable method for distributing the RBF centers as regions with more samples will also have more spatial resolution.

There are, of course, other dictionaries that may prove more effective in other circumstances.
For example, basis functions defined in polar coordinates are useful when limit cycles or other periodic orbits are present as they mimic the form of the Koopman eigenfunctions for simple limit cycles~\cite{Bagheri2013}.
How to choose the best set of trial functions is an important, yet open, question; fortunately, the EDMD method often produces useful results even with the relatively naive choices of trial functions presented in this section.

\section{Deterministic Data and the Koopman Eigenfunctions}
\label{sec:deterministic}

Most applications of DMD assume that the data sets were generated by a deterministic dynamical system.
In Sec.~\ref{sec:dmd}, we showed that EDMD produces an approximation of the Koopman eigenfunctions, eigenvalues, and modes with large amounts of data.
In this section, we demonstrate that EDMD can produce accurate approximations of the Koopman eigenfunctions, eigenvalues, and modes with limited amounts of data by applying the method to two illustrative examples.
The first is a discrete time linear system, one where the eigenfunctions, eigenvalues, and modes are known analytically, and serves as a test case for the method.
The second is the unforced Duffing oscillator.
Our goal there is to demonstrate that the approximate Koopman eigenfunctions obtained via EDMD have the potential to serve as a data driven parameterization of a system with multiple basins of attraction.

\subsection{A Linear Example}
\label{subsec:linear}

\subsubsection{The Governing Equation, Data, and Analytically Obtained Eigenfunctions}

One system where the Koopman eigenfunctions are known analytically is a simple LTI system of the form
\begin{equation}
\vec x(n+1) = \mat{J} \vec x(n),
\end{equation}
with $\vec x(n)\in\mathbb{R}^N$ and $\mat{J}\in\mathbb{R}^{N\times N}$.
It is clear that an eigendecomposition  yields complete information about the underlying system provided $\mat{J}$ has a complete set of eigenvectors.
Because the underlying dynamics are linear, it should not be surprising that the Koopman approach contains {\em the eigendecomposition}.

To show this, note that the action of the Koopman operator for this problem is 
\begin{equation}
\koopman\phi(\vec x) = \phi(\mat{J}\vec x),
\label{eq:koopman-linear}
\end{equation}
where $\phi\in\observables$.  
Assuming $\mat{J}$ has a complete set of eigenvectors, it will have $N$ left eigenvectors, $\vec w_i$, that satisfy $\vec w_i^*\mat{J} = \mu_j\vec w_i^*$, and where the $i$-th eigenvector is associated with the eigenvalue $\mu_i$.   
Then, the function
\begin{equation}
 \varphi_{n_1, n_2, \ldots, n_N}(\vec x) = \prod_{i=1}^N (\vec w_i^*\vec x)^{n_i}
\label{eq:koopman-eigenfunction-linear}
\end{equation}
is an eigenfunction of the Koopman operator with the eigenvalue $\prod_{i=1}^N \mu_i^{n_i}$ for $n_i\in\mathbb{N}$.  
This is a well known result (see, e.g., Ref.~\cite{Budivsic2012}), so we show the proof of this only for the sake of completeness. 
We proceed directly:
\begin{equation*}
\begin{aligned}
 \koopman \varphi_{n_1, n_2, \ldots, n_N}(\vec x) &= \prod_{i=1}^N (\vec w_i^*\mat{J}\vec x)^{n_i}
= \prod_{i=1}^N \mu_i^{n_i}(\vec w_i^*\vec x)^{n_i} \\
& = \left(\prod_{i=1}^N \mu_i^{n_i}\right)\varphi_{n_1, n_2, \ldots, n_N}(\vec x),
\end{aligned}
\end{equation*}
by making use of \eqref{eq:koopman-linear} and the definition of $\vec w_i$ as a left eigenvector. 
Then, the representation of the full state observable in terms of the Koopman modes and eigenfunctions (i.e., \eqref{eq:koopman-modes}) is 
\begin{equation}
\vec x = \sum_{i=1}^N \vec v_i \left(\vec w_i^*\vec x\right), 
\label{eq:linear-koopman-modes}
\end{equation}
where the $i$-th Koopman mode, $\vec v_i$, is the $i$-th eigenvector of $\mat{J}$ suitably scaled so that $\vec w_i^*\vec v_i = 1$.
This is identical to writing $\vec x$ in terms of the eigenvectors of $\mat{J}$; inner products with the left eigenvectors determine the component in each direction, and the (right) eigenvectors  allow the full state to be reconstructed.

\begin{figure}
\centering
\includegraphics[width=\textwidth]{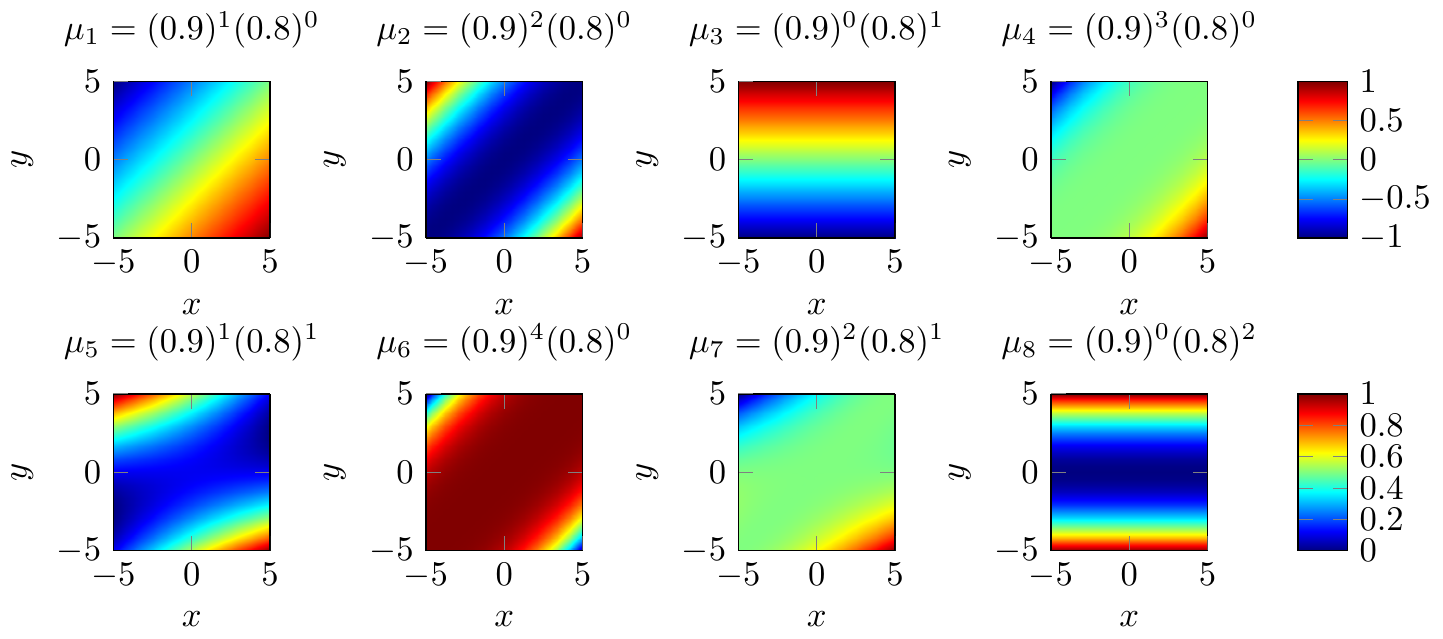}
\caption{Pseudocolor plots of the first 8 Koopman eigenfunctions plotted using their analytical expressions, \eqref{eq:linear-eigenfunctions}, along with the associated eigenvalue.
Note that the eigenfunctions were scaled so that $\|\varphi_i\|_\infty=1$ on the domain shown.
}
\label{fig:koopman-linear-eigenfunctions}
\end{figure}

As a concrete example, consider 
\begin{equation}
\vec x(n+1) = 
\begin{bmatrix}
        0.9 & -0.1 \\
        0.0 &  0.8
\end{bmatrix}\vec x(n) = \mat{J}\vec x(n),
\label{eq:linear}
\end{equation}
where $\vec x_{n} = [x_n,y_n]$.
From \eqref{eq:koopman-eigenfunction-linear}, the Koopman eigenfunctions and eigenvalues are 
\begin{equation}
\varphi_{ij}(x, y) = \left(\frac{x - y}{\sqrt 2}\right)^iy^j, \qquad 
\lambda_{ij} = (0.9)^i(0.8)^j,
\label{eq:linear-eigenfunctions}
\end{equation}
for $i,j\in\mathbb{Z}$.  
Figure~\ref{fig:koopman-linear-eigenfunctions} shows the first 8 nontrivial eigenfunctions sorted by their associated eigenvalue. 
The 0-th eigenfunction, $\varphi_{00}(\vec x) = 1$ with $\mu_{00} = 1$, was omitted because it is always an eigenfunction and will be recovered by EDMD if $\psi=1$ is included as a dictionary element.

To apply the EDMD procedure, one needs both data and a dictionary of observables.  
The data in $\mat{X}$ consists of 100 normally distributed initial conditions, $\vec x_i$, and their images, $\vec y_i = \mat{J}\vec x_i$, which we aggregate in the matrix $\mat{Y}$, i.e., $\mat{X},\mat{Y}\in\mathbb{R}^{2\times 100}$.
The dictionary, $\dictionary$, is chosen to contain the direct sum of Hermite polynomials in $x$ and $y$ that include up to the 5th order terms in $x$ and $y$, i.e., 
\begin{equation}
\begin{aligned}
\mathcal{D} &= \{\psi_0, \psi_1, \psi_2, \ldots\} \\
& = \small\begin{matrix}
\{H_0(x)H_0(y), & H_1(x)H_0(y), & H_2(x)H_0(y), & H_3(x)H_0(y), &H_4(x)H_0(y), \\
H_0(x)H_1(y), & H_1(x)H_1(y), & H_2(x)H_1(y), & H_3(x)H_1(y), &H_4(x)H_1(y), \\
H_0(x)H_2(y), & H_1(x)H_2(y), & H_2(x)H_2(y), & H_3(x)H_2(y), &H_4(x)H_2(y), \\
H_0(x)H_3(y), & H_1(x)H_3(y), & H_2(x)H_3(y), & H_3(x)H_3(y), &H_4(x)H_3(y), \\
H_0(x)H_4(y), & H_1(x)H_4(y), & H_2(x)H_4(y), & H_3(x)H_4(y), &H_4(x)H_4(y)\},
\end{matrix}
\end{aligned}
\label{eq:basis-ordering}
\end{equation}
where $\dictionary$ has been written this way to make the ordering of the dictionary elements apparent.
The Hermite polynomials were chosen because they are an appropriate basis for Cauchy problems, and orthogonal with respect to the weight function, $\rho(\vec x) = e^{-\|\vec x\|^2}$, that is implicit in the normally distributed sampling strategy used here.

\subsubsection{Results}

\begin{figure}
\centering
\includegraphics[width=\textwidth]{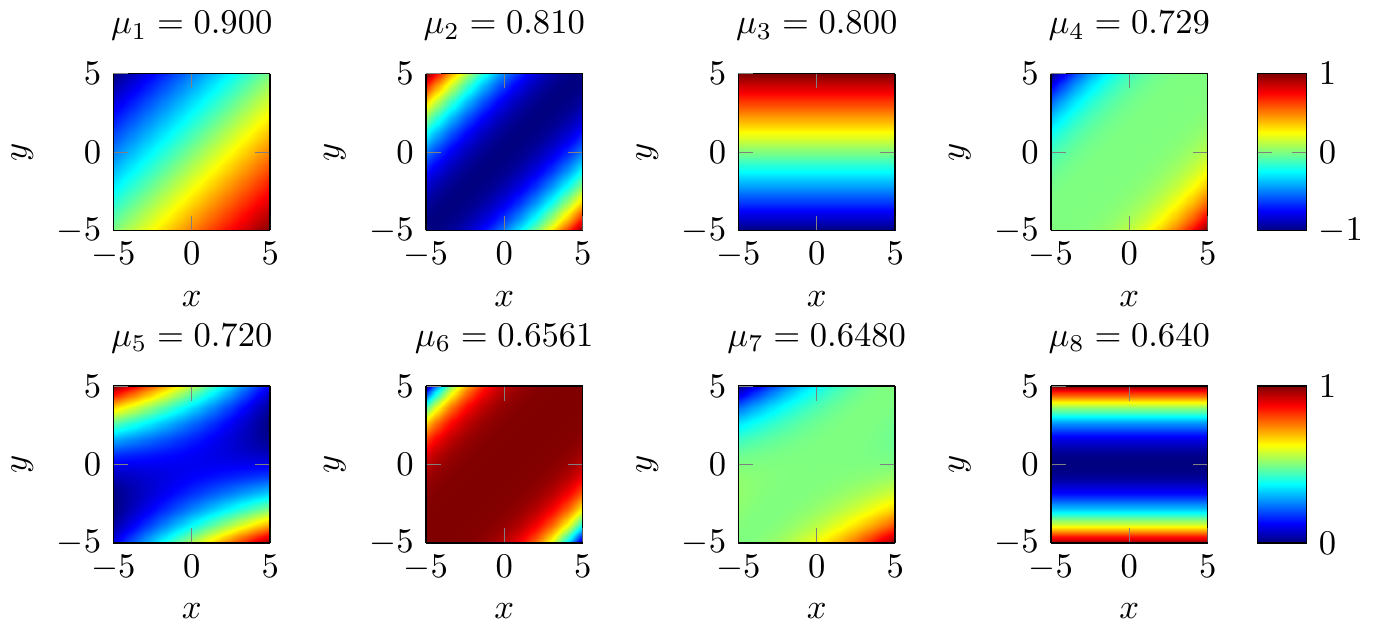}
\caption{
Pseudocolor plots of the first eight nontrivial eigenfunctions of the Koopman operator computed using the EDMD procedure.
The eigenfunctions obtained using EDMD were scaled such that $\|\varphi_i\|_\infty=1$ on the domain shown.
With this scaling, there is excellent agreement between these results and those presented in Fig.~\ref{fig:koopman-linear-eigenfunctions}.
}
\label{fig:edmd-eigenfunctions}
\end{figure}

Figure~\ref{fig:edmd-eigenfunctions} shows the same 8 eigenfunctions computed using the EDMD method.
Overall, there is (as expected) excellent quantitative agreement, both in the eigenvalues and the eigenfunctions, with the analytical results presented in Fig.~\ref{fig:koopman-linear-eigenfunctions}.
On the domain shown, the eigenvalues are accurate to 10 digits, and the maximum pointwise difference between the true and computed eigenfunction is $10^{-6}$.  
In this problem, standard DMD also generates highly accurate approximations of $\varphi_1$ and $\varphi_3$ and their associated eigenvalues, but {\em will not produce any of the other eigenfunctions}; the standard choice of the dictionary only contains linear terms and, therefore, cannot reproduce eigenfunctions with constant terms or any nonlinear terms.
As a result, expanding the basis allows EDMD to capture more of the Koopman eigenfunctions than standard DMD could.
These additional eigenfunctions are not necessary for an LTI system, but are in principle needed in nonlinear settings.

\begin{figure}
\centering
\includegraphics[width=\textwidth]{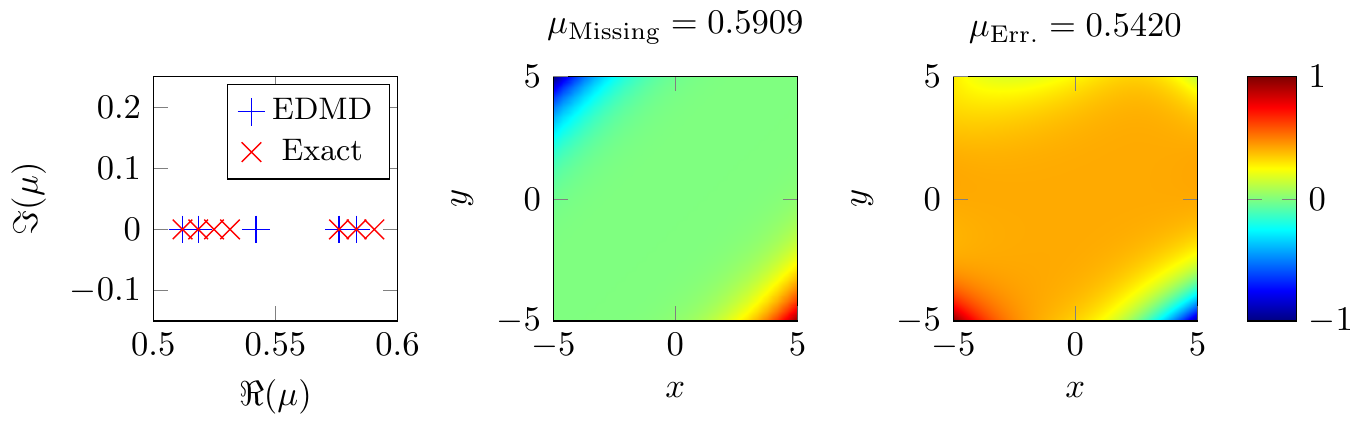}
\caption{
A subset of the spectrum of the Koopman operator and the EDMD computed approximation. As shown here, there are clearly errors in the spectrum further away from the unit circle (which is not contained in the plotting region).
The center plot shows an example of a ``missing'' eigenfunction that is not captured by EDMD; this eigenfunction is $(x-y)^5\not\in\observables_{\dictionary}$, and cannot be represented with our basis functions.
The right plot shows an example of an ``erroneous'' eigenfunction that appears because $\observables_{\dictionary}$ is not an invariant subspace of the Koopman operator. 
}
\label{fig:linear-eigenvalues}
\end{figure}

This level of accuracy is in large part because the first nine eigenfunctions are in $\observables_{\dictionary}$, the subspace of observables spanned by $\dictionary$.
When this is not the case, the result is either a {\em missing or erroneous} eigenfunctions like the examples shown in Fig.~\ref{fig:linear-eigenvalues}.
The 9-th eigenfunction, $(x-y)^5$ with $\mu = 0.9^5 = 0.59049$, is not captured by EDMD {\em with the dictionary chosen here} because it lacks the needed 5-th order monomials in $x$ and $y$, which is similar to how DMD skips the 2nd Koopman eigenfunction due to a lack of quadratic terms.

The erroneous eigenfunction appears because $\observables_{\dictionary}$ is not invariant with respect to the action of the Koopman operator. 
In particular, $\phi$ contains the term $yx^4$ whose image $\koopman(yx^4)\not\in\observables_{\dictionary}$ because $x^5,y^5\not\in\observables_{\dictionary}$. 
In most applications, there are small components of the eigenfunction that cannot be represented in the dictionary chosen, which results in errors in the eigenfunction such as the one seen here. 
Even in the limit of infinite data, we would compute the eigenfunctions of $\mathcal{P}_{\observables_\dictionary}\koopman$, where $\mathcal{P}_{\observables_\dictionary}$ is the projection onto $\observables_{\dictionary}$, rather than the eigenfunction of $\koopman$.
To see that this not a legitimate eigenfunction, we added $H_5(x)$ and $H_5(y)$ to $\dictionary$, which removes this erroneous eigenfunction. 

Finally, we compute the Koopman modes for the full state observable.
Using the ordering of the dictionary elements given in \eqref{eq:basis-ordering}, the weight matrix, $\mat{B}$ in \eqref{eq:vector-valued-basis}, needed to compute the Koopman modes for the full state observable, $\vec g(\vec x)^T = [x, y]$ is 
\begin{equation}
\mat{B} = 
\begin{bmatrix}
0 & 0 \\
1 & 0 \\
0 & 0 \\
0 & 0 \\
0 & 0 \\
0 & 1 \\
0 & 0 \\
\vdots & \vdots
\end{bmatrix}.
\end{equation}  
The Koopman modes associated with $\mu_1 = 0.9$ is $\vec v_1 = [0, -\sqrt{2}]^T$, while the Koopman mode associated with $\mu_3 = 0.8$ is $\vec v_3 = [-1, -1]$; again, these are the eigenvectors of $\mat J$.
The contribution of the other {\em numerically computed} eigenfunctions in reconstructing the full state observable is negligible (i.e., $\|\vec v_k\|\approx 10^{-11}$ for $k\neq 1,3$), so the Koopman/EDMD analysis is  an eigenvalue/eigenvector decomposition once numerical errors are taken into consideration. 
Although EDMD reveals a richer set of Koopman eigenfunctions that are analytically known to exist, their associated Koopman modes are zero and, hence, they can be neglected.
Our goal in presenting this example is not to demonstrate any new phenomenon, but rather to demonstrate that there is good quantitative agreement between the analytically obtained Koopman modes, eigenvalues, and eigenfunctions and the approximations produced by EDMD. 
Furthermore, it allowed us to highlight the types of errors that appear when  $\observables_{\dictionary}$ is not an invariant subspace of $\koopman$, which results in erroneous eigenfunctions, or when the dictionary is missing elements, which results in missing eigenfunctions.

\subsection{The Duffing Oscillator}
\label{subsec:duffing}

In this section, we will compute the Koopman eigenfunctions for the unforced Duffing Oscillator, which for the parameter regime of interest here, has two stable spirals and a saddle point whose stable manifold defines the boundary between the basins of attraction.
Following  Ref.~\cite{Mauroy2013} and the references contained therein, the eigenvalues of the linearizations about the fixed points in the system are known to be a subset of the Koopman eigenvalues, and for each  stable spiral, the magnitude and phase of the associated Koopman eigenfunction parameterizes the relevant basin of attraction.
Additionally, because basins of attraction are forward invariant sets, there will be two eigenfunctions with $\mu = 0$, each of which is supported on one of the two basins of attraction in this system (or, equivalently, there will be a trivial eigenfunction and another eigenfunction with $\mu=0$ whose level sets denote the basins of attraction). 
Ultimately, we are not interested in recovering highly accurate eigenfunctions in this example. 
Instead, we will demonstrate that the eigenfunctions computed by EDMD are {\em accurate enough} that they can be used to identify and parameterize the basins of attraction that are present in this problem for the region of interest.

The governing equations for the unforced Duffing Oscillator are 
\begin{equation}
\ddot x = -\delta \dot x - x(\beta + \alpha x^2),
\end{equation}
which we will study using the parameters $\delta = 0.5$, $\beta = -1$, and $\alpha = 1$.
In this regime, there are two stable spirals at $x = \pm 1$ with $\dot x = 0$, and a saddle at $x,\dot x = 0$, so almost every initial condition (except for those on the stable manifold of the saddle) is drawn to either of the spirals. 
In what follows, the data consist of  $10^3$ trajectories with 11 samples each with  
a sampling interval of $\Delta t = 0.25$ (i.e., $\mat{X},\mat{Y}\in\mathbb{R}^{2\times 10^4}$), and initial conditions uniformly distributed over $x,\dot x\in[-2, 2]$. 
With this sampling rate and initialization scheme, many trajectories will {\em approach} the stable spirals, but few will have (to numerical precision) reached the fixed points.
As a result, the basins of attraction cannot be determined by observing the last snapshot in a given trajectory.
Instead, EDMD will be used to ``stitch'' together this ensemble of trajectories to form a single coherent picture.

However, because there are multiple basins of attraction, the leading eigenfunctions will be discontinuous~\cite{Mauroy2013}, and supported only on the appropriate basin of attraction.
In principle, our computation could be done ``all at once'' using a single $\dictionary$ and applying EDMD to the complete data set.
To enforce the compactly supported nature of the eigenfunctions regardless of which dictionary we use, we will proceed in a two-tiered fashion. 
First, the basins of attraction will be identified using all of the data and a dictionary with support everywhere we have data. 
Once we have identified these basins, both state space and the data will be partitioned into subdomains based on the numerically identified basins. 
The EDMD procedure will then be run on each subdomain and the corresponding partitioned data set individually.

\subsubsection{Locating Basins of Attraction}

\begin{figure}
\centering
\includegraphics[width=\textwidth]{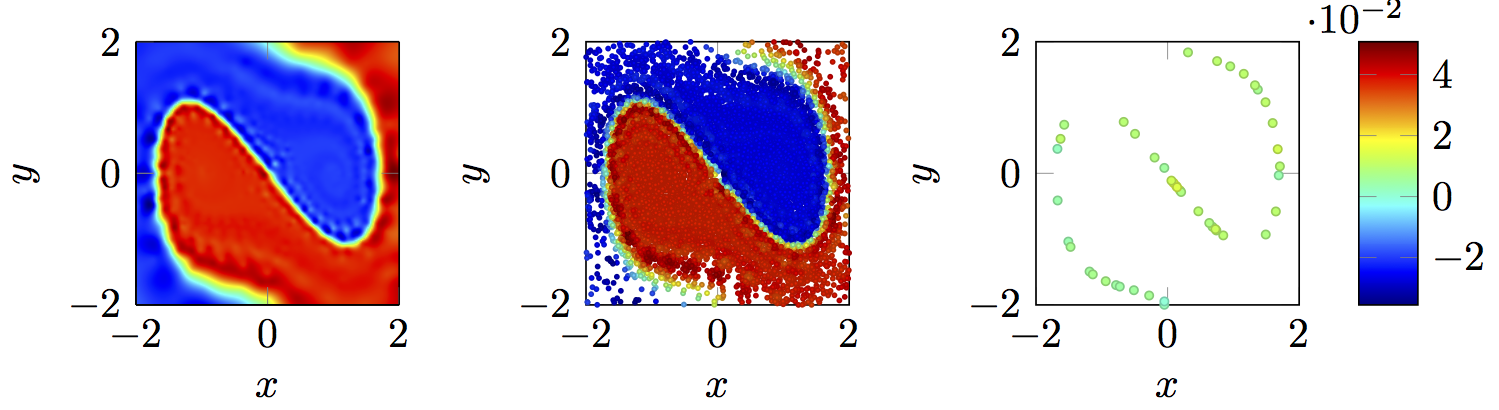}
\caption{
(left) A plot of the first nontrivial eigenfunction generated with an eigenvalue of $\lambda = 0.048$ obtained from $10^3$ randomly initialized trajectories each consisting of 10 steps taken with $\Delta t = 0.25$.  
This eigenfunction should be constant in each basin of attraction and have $\lambda = 0$, so EDMD is generating an {\em approximation of the eigenfunction} rather than a fully converged solution.
(center) The first nontrivial eigenfunction evaluated at the data points.  
(right) A plot of the mis-classified points (less than $0.5$\% of the data), each of these points are (to the eye) on the boundary between the invariant sets.
}
\label{fig:duffing-basin}
\end{figure}

Figure~\ref{fig:duffing-basin} highlights the first step: partitioning of state space into basins of attraction.
We used a dictionary consisting of the constant function and  1000 radial basis functions (RBFs) (the thin plate splines described in Sec.~\ref{sec:dictionary}), where $k$-means clustering~\cite{Bishop2006} on the full data set was used to choose the RBF centers.
RBFs were chosen here because of the  geometry of the computational domain; indeed, RBFs are often a fundamental component of ``mesh-free'' methods that avoid the nontrivial task of generating a computational mesh~\cite{Liu2010}.

The  leading (continuous time) eigenvalue is $\lambda_0 = -10^{-14}$ which corresponds to the constant function.
The second eigenfunction, shown in the leftmost image of Fig~\ref{fig:duffing-basin}, has $\lambda_1 = -10^{-3}$, which should  be considered an approximation of zero.
The discrepancy between the numerically computed eigenfunction and the theoretical one is due to the choice of the dictionary.
The analytical eigenfunction possesses a discontinuity on the edge of the basin of attraction (i.e., the stable manifold of the saddle point at the origin), but discontinuous functions are not in the space spanned by RBFs.
Therefore,  the numerically computed approximation ``blurs'' this edge as shown in Fig~\ref{fig:duffing-basin}. 

The scatter plot in the center of Fig.~\ref{fig:duffing-basin} shows the data points colored by the 1st nontrivial eigenfunction.
There is good qualitative agreement between the numerically computed basin of attraction and the actual basin.
By computing the mean value of $\varphi_1$ on the data and using that value as the threshold that determines which basin of attraction a point belongs to, the EDMD approach mis-classifies only 46 of the $10^4$ data points, resulting in an error of only 0.5\% as shown by the rightmost plot.  
As a result, the leading eigenfunctions computed by EDMD are sufficiently accurate to produce a meaningful partition of the data.

\subsubsection{Parameterizing a Basin of Attraction}

\begin{figure}
\centering 
\includegraphics[width=\textwidth]{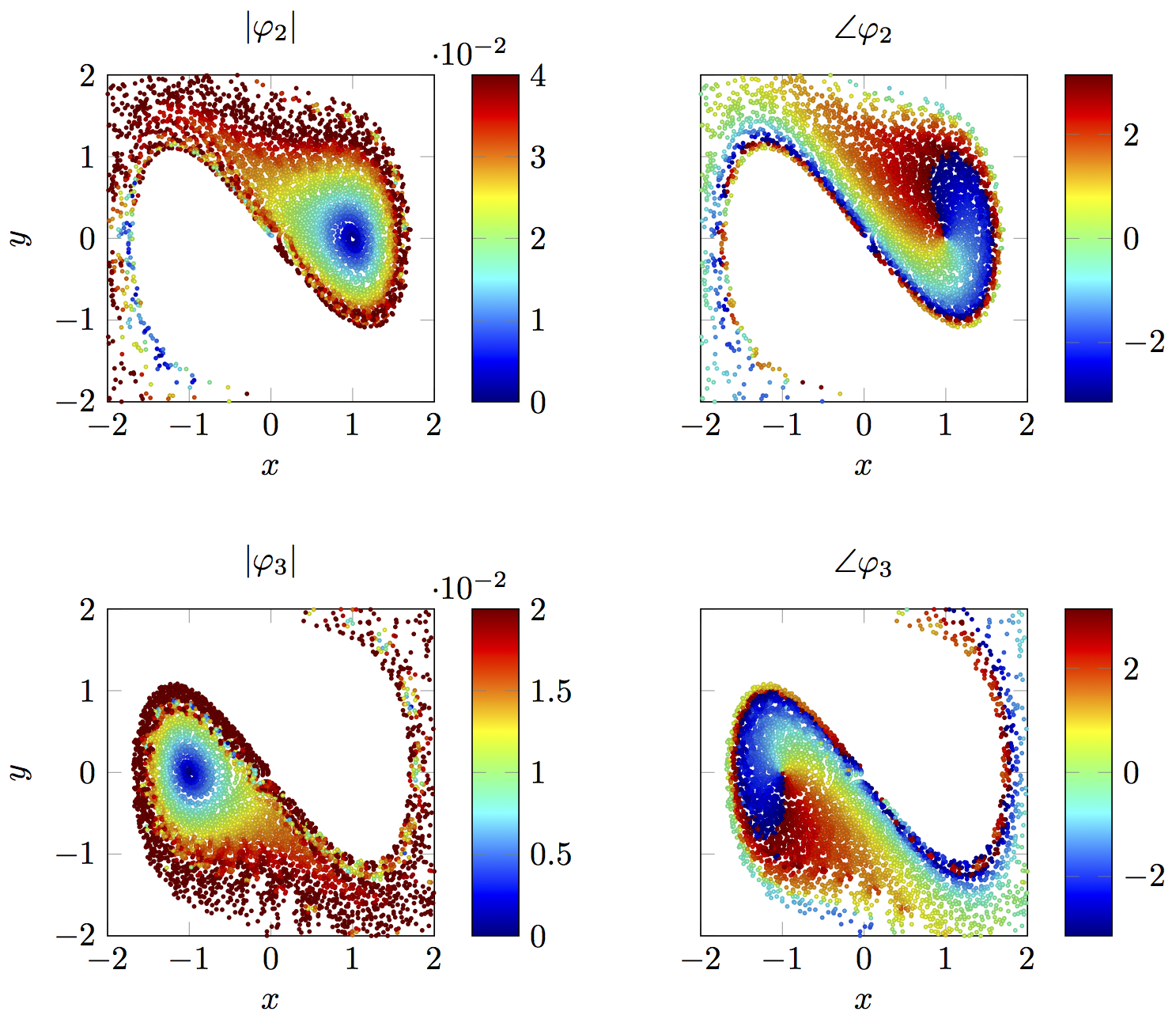}
\caption{(top) The amplitude and phase of the Koopman eigenfunction with $\lambda = -0.237 + 1.387\imath$ (analytically, $-0.25 + 1.3919\imath$) for the stable spiral at $(1,0)$.
(bottom) The same pair of plots for the spiral at $(-1, 0)$. 
Although there are errors near the ``edge'' of the basin of attraction, the amplitude and phase of this eigenfunction can serve as a polar coordinate system for the corresponding basin of attraction. 
 }
\label{fig:duffing-basin-parameterization}
\end{figure}

Now that the basins of attraction have been identified, the next task is to develop a coordinate system or parameterization of the individual basins. 
To do so, we will use the eigenfunctions associated with the eigenvalues of the system linearization about the corresponding fixed point.
Because these fixed points are spirals, this parameterization can be realized using the amplitude and phase of one member of the complex-conjugate pair of eigenfunctions. 
To approximate these eigenfunctions, we first partition our data into two sets as mentioned above using the leading Koopman eigenfunctions (including the mis-classified data points). 
On each subset of the data, the $k$-means procedure was run again to select a new set of 1000 RBF centers, and this ``adjusted'' basis along with the constant function comprised the $\dictionary$ used by EDMD. 
Figure~\ref{fig:duffing-basin-parameterization} shows the amplitude and phase of the eigenfunction with eigenvalue closest to $-0.25 + 1.387\imath$ computed using the data in each basin of attraction. 
The computed eigenvalues agree favorably with the analytically obtained eigenvalue; the basin of the spiral at $(1,0)$ has the eigenvalue $-0.237 + 1.387\imath$, and the basin of the spiral at $(-1, 0)$ has the eigenvalue $-0.24 + 1.35\imath$.  

Figure~\ref{fig:duffing-basin-parameterization} demonstrates that the amplitude and phase of a Koopman eigenfunction forms something analogous to an ``action--angle'' parameterization of the basin of attraction.
Due to the nonlinearity in the Duffing oscillator, this parameterization is more complicated than an appropriately shifted polar coordinate system, and is, therefore, not the parameterization that would be generated by linearization about either $(\pm 1, 0)$.
The level sets of the amplitude of this eigenfunction are the so-called ``isostables''~\cite{Mauroy2013}.
One feature predicted in that manuscript is that the 0-level set of the isostable is the fixed point in the basin of attraction; this feature is reflected in Fig.~\ref{fig:duffing-basin-parameterization} by the blue region, which corresponds to small values of the eigenfunction that are near the fixed points at $(\pm 1, 0)$. 
Additionally, a singularity in the phase can be observed there. 
The EDMD approach produces noticeable numerical errors near the edges of the basin.
These errors can be due to a lack of data, or to the singularities in the eigenfunctions that can occur at unstable fixed points~\cite{Mauroy2013a}.

In this section, we applied the EDMD procedure to deterministic systems and showed that it produces an approximation of the Koopman operator.
With a sensible choice of data and $\dictionary$, we showed that EDMD generates a quantitatively accurate approximation of the Koopman eigenvalues, eigenfunctions, and modes for the linear example. 
In the second example, we used the Koopman eigenfunctions to identify and parameterize the basins of attraction of the Duffing Oscillator. 
Although the EDMD approximation of the eigenfunctions could be made more accurate with more data,  it is still accurate enough to serve as an effective parameterization. 
As a result, the EDMD method can be useful outside of the large data limit, and should be considered an enabling technology for data driven approximations of the Koopman eigenvalues, eigenfunction, and modes.

\section{Stochastic Data and the Kolmogorov Backward Equation}
\label{sec:stochastic-ck}

The EDMD approach is entirely data driven, and will produce an output regardless of the nature of the data given to it. 
However, if the results of EDMD are to be meaningful, then certain assumptions must be made about the dynamical system that produced the data used.
In the previous section, it was assumed the data were generated by a deterministic dynamical system; as a result, EDMD produced approximations of the tuples of Koopman eigenfunctions, eigenvalues, and modes. 

Another interesting case to consider is if the underlying dynamical system is a Markov process, such as a stochastic differential equation (SDE).
For such systems, the evolution of an observable is governed by the Kolmogorov backward (KB) equation~\cite{Bagherisubmitted}, whose ``right hand side'' has been called the ``stochastic Koopman operator'' (SKO)~\cite{Mezic2005}.
In this section, we will show that EDMD produces approximations of the eigenfunctions, eigenvalues, and modes of the SKO if the underlying dynamical system happens to be a Markov process.

To accomplish this, we will prove that the EDMD method converges to a Galerkin method in the large data limit.
After that, we will demonstrate its accuracy with finite amounts of data by applying it to the model problem of a 1D SDE with a double well potential, where the SKO eigenfunctions can be computed using standard numerical methods.

Another proposed application of the Koopman operator is for the purposes of {\em model reduction}, which as been explored  in Ref.~\cite{Froyland2013}.
 Model reduction based on the Koopman eigenfunctions is equally applicable in both deterministic and stochastic settings, but we choose to present it for stochastic systems to highlight the similarities between EDMD and {\em manifold learning techniques} such as diffusion maps~\cite{Nadler2005, Coifman2006}.
In particular, we apply EDMD to an SDE defined on a ``Swiss Roll,'' which is a nonlinear manifold often used to test manifold learning methods~\cite{Lee2007}.
The purpose of this example is twofold: first, we show that a data driven parameterization of the Swiss Roll can be obtained using EDMD, and second, we show that this parameterization will preferentially capture ``slow'' dynamics on that manifold before the ``fast'' dynamics when the noise is made anisotropic.

\subsection{EDMD with Stochastic Data}
\label{subsec:stochastic-edmd}
For a {\em discrete time Markov process}, 
\begin{equation*}
        \vec{x} \mapsto\vec{F}(\vec{x};\vec{\omega}),
\end{equation*}
the SKO~\cite{Mezic2005} is defined as
\begin{equation}
(\tilde\koopman\psi)(\vec{x})=\expect[\psi(\vec{F}(\vec{x};\vec{\omega}))],\label{eq:koopman-stochastic}
\end{equation}
where $\omega\in\Omega_s$ is an element in the probability space associated with the stochastic dynamics ($\Omega_s$),  $\expect$ denotes the expected value over that space, and $\psi\in\observables$ is a scalar observable. 
The SKO  ~\cite{Mezic2005} takes an observable of the full system state and returns the conditional expectation of the observable ``one timestep in the future.''
Note that this definition is compatible with the deterministic Koopman operator because $\expect[\psi(\vec{F}(\vec{x}))]=\psi(\vec{F}(\vec{x}))$ if $\vec F$ is deterministic.

As with the deterministic case, we assume the snapshots in $\mat{X}\in\mathbb{R}^{N\times M}$ were generated by randomly placing initial conditions on $\manifold$ with the  density of $\rho(\vec x)$ and that $M$ is sufficiently large. 
Once again, $\rho$ does not need to be an invariant measure of the underlying dynamical system; it is simply the sampling density of the data.
Due to the stochastic nature of the system, there are two probability spaces involved: one related to the samples in $\mat X$ and
another for the stochastic dynamics.
Because our system has ``process'' rather than ``measurement'' noise, the $\vec x_i$ are known exactly, and the interpretation of the Gram matrix, $\mat{G}$, remains unchanged.
Therefore,
\[
\lim_{M\to\infty} {\mat{G}}_ {i,j}=\int_{\manifold}\psi_{i}^*(\vec{x})\psi_{j}(\vec{x})\rho(\vec{x})d\vec{x}=\left\langle \psi_{i},\psi_{j}\right\rangle_{\rho}, 
\]
by the law of large numbers when $M$ is large enough.
This is identical to the deterministic case.
However, the definition of  $\mat A$ will change.
Assuming that the choice of $\vec \omega$ and $\vec x$ are independent, 
\begin{equation*}
\begin{aligned}
\lim_{M\to\infty}\mat A_{i,j}&=\expect[\psi_{i}^*(\tilde\koopman\psi_{j})] 
=\int_{\manifold\times\Omega_{s}}\psi_{i}^*(\vec{x})\psi_{j}(\vec{F}(\vec{x},\vec{\omega}))\rho(\vec{x})_{}\; d\vec{x}d\vec{\omega} \\ &=\int_{\manifold}\psi_{i}^*(\vec{x})\expect[\psi_{j}(\vec{F}(\vec{x},\vec{\omega}))]\rho(\vec{x})\, d\vec{x}
=\left\langle \psi_{i},\tilde\koopman\psi_{j}\right\rangle_{\rho},
\end{aligned}.
\end{equation*}
The elements of $\mat A$ now contain a {\em second} integral over the probability space that pertains to the stochastic dynamics, which produces the  expectation of the observable in the expression above.

The accuracy of the resulting method will depend on the dictionary, $\dictionary$, the manifold on which the dynamical system is defined, the data, and the dynamics used to generate it.  
One interesting special case is if the basis functions are indicator functions supported on  ``boxes.''
When this is the case, EDMD is equivalent to the widely used Ulam Galerkin method~\cite{Dellnitz2001,Bollt2013}.
This equivalence is lost for other choices of $\dictionary$ and $\rho$, but as we will demonstrate in the subsequent sections, EDMD can produce accurate approximations of the eigenfunctions for many other choices of these quantities.

The ``stochastic Koopman modes'' can then be computed using \eqref{eq:compute-modes}, but they too must be reinterpreted as the weights needed to reconstruct the {\em expected value} of the full state observable using the eigenfunctions of the SKO. 
Due to the stochastic nature of the dynamics, the Koopman modes can no longer exactly specify the state of the system.
However, they can be used as approximations of the Koopman modes that would be obtained in the ``noise free'' limit when some appropriate restrictions are placed on the nature of the noise and the underlying dynamical system.
Indeed, these are the modes we are truly computing when we apply DMD or EDMD to experimental data, which by its very nature, contains some noise.

\subsection{A Stochastic Differential Equation with a Double Well Potential }
\label{subsec:double_well}

In this section, we will show that the EDMD procedure is capable of accurately approximating the eigenfunctions of the stochastic Koopman operator by applying it to an SDE with a double well potential. 
Although we do not have analytical solutions for the eigenfunctions, the problem is simple enough that we can accurately compute them using
standard numerical methods.

\subsubsection{The Double Well Problem and Data}
\begin{figure}
\centering
\includegraphics{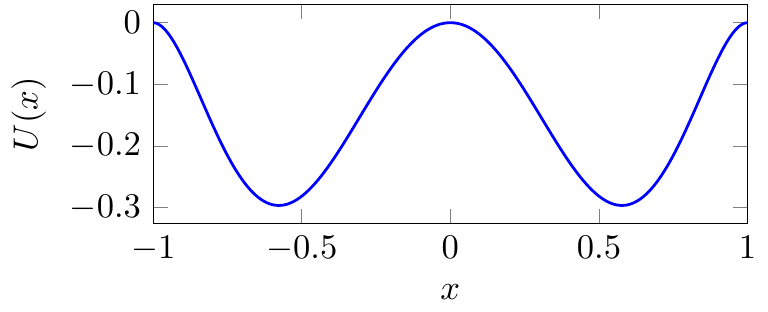}
\caption{The double well potential $U(x) = -2(x^2 - 1)^2x^2$ }
\label{fig:wells}
\end{figure}

First, consider an SDE with a double well potential.
Let the governing equations for this system be 
\begin{equation}
d x = -\nabla U(x) dt +  \sigma d W_t,
\label{eq:sde}
\end{equation}
where $x$ is the state, $-\nabla U(x)$ the drift, and $ \sigma$ is the (constant) the diffusion coefficient.
Furthermore, no flux boundary conditions are imposed at $x = \pm 1$. 
For this problem, we let $U(x) = -2 (x^2 - 1)^2x^2$ as shown in Fig.~\ref{fig:wells}.

The Fokker-Planck equation associated with this SDE is 
\begin{equation}
\frac{\partial \rho(x,t)}{\partial t} 
         = - \frac{\partial}{\partial x}\left(-\frac{\partial U}{\partial x} \rho(x,t)\right)
          + \frac{\sigma^2}{2}\frac{\partial^2\rho(x,t)}{\partial x^2}
          = \mathcal{P}\rho,
\end{equation}
where $\rho$ is a probability density with $\partial_x\rho(x,t)\big|_{x=\pm 1} = 0$ due to the no-flux boundary conditions we impose, and $\mathcal{P}$ is the Perron-Frobenius operator. 
The adjoint of the Perron-Frobenius operator determines the Kolmogorov backward equation, and thus defines the stochastic Koopman operator, $\tilde\koopman = \mathcal{P}^\dagger$. 
For this example, 
\begin{equation}
\tilde\koopman\phi =  -\frac{\partial U}{\partial x}\frac{\partial\phi}{\partial x} + 
        \frac{\sigma^2}{2}\frac{\partial^2 \phi}{\partial x^2} 
        \label{eq:stochastic-koopman-wells}
\end{equation}
with Neumann boundary conditions, $\partial_x\phi\big|_{x=\pm 1} = 0$.
To directly approximate the Koopman eigenfunctions,
\eqref{eq:stochastic-koopman-wells} is discretized in space using a second order finite difference scheme with 1024 interior points. 
The eigenvalues and eigenvectors of the resulting finite dimensional approximation of the Koopman operator will be used to validate  the EDMD computations.

The data are $10^6$ initial points on $x\in [-1, 1]$ drawn from a uniform distribution, which constitute  $\mat{X}$, and their positions after $\Delta t = 0.1$, which constitute $\mat{Y}$.
The evolution of each initial condition was accomplished through $10^2$ steps of the Euler--Maruyama method~\cite{Higham2001, Kloeden1992} with a timestep of $10^{-3}$ using the double well potential in Fig.~\ref{fig:wells}. 
The dictionary chosen is a discontinuous spectral element basis that
splits $x\in[-1,1]$ into four equally sized subdomains with up to tenth
order Legendre polynomials on each subdomain (see Sec.~\ref{sec:dictionary}) for a total of forty degrees of freedom. 

\subsubsection{Recovering the Koopman Eigenfunctions and Eigenvalues}

\begin{figure}
\centering
\includegraphics[width=\textwidth]{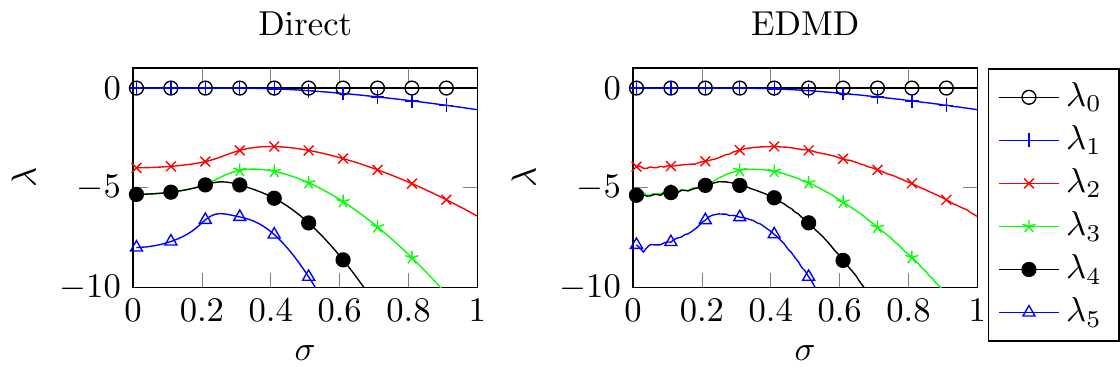}
\caption{(left) The first six eigenvalues of the stochastic Koopman operator obtained with a finite difference discretization of $\tilde\koopman$.
(right) The first six eigenvalues obtained with the EDMD approach. 
In both plots, a marker is placed every 10-th data point.
While there is good quantitative agreement between the true eigenvalues and those obtained with EDMD, some small ``noise'' due to quadrature errors does appear in the right plot as $\sigma\to 0$.
}
\label{fig:double-well-eigenvalues}
\end{figure}

Because the Koopman operator is infinite dimensional, we will clearly be unable to approximate all of the tuples. 
Instead, we focus on the leading (i.e., most slowly decaying) tuples, which  govern the long--term dynamics of the underlying system.
In this example, we seek to demonstrate that our approximation is: (a) quantitatively accurate, and (b) valid over a range of coefficients, $\sigma$, and not solely in the small (or large) noise limits.

Figure~\ref{fig:double-well-eigenvalues} shows the first six eigenvalues obtained using a finite difference discretization of the Koopman operator and the EDMD approximation as a function of $\sigma$. 
In this problem, $\varphi_0=1$, is always an eigenfunction of the Koopman operator with $\lambda_0 = 0$ for all values of $\sigma$.
Because  $\dictionary$ contains the constant function, it should be no surprise the EDMD method is able to identify it as an eigenfunction.
Though trivial, the existence of this eigenfunction is a ``sanity check'' for the method.

\begin{figure}
\centering
\includegraphics[width=\textwidth]{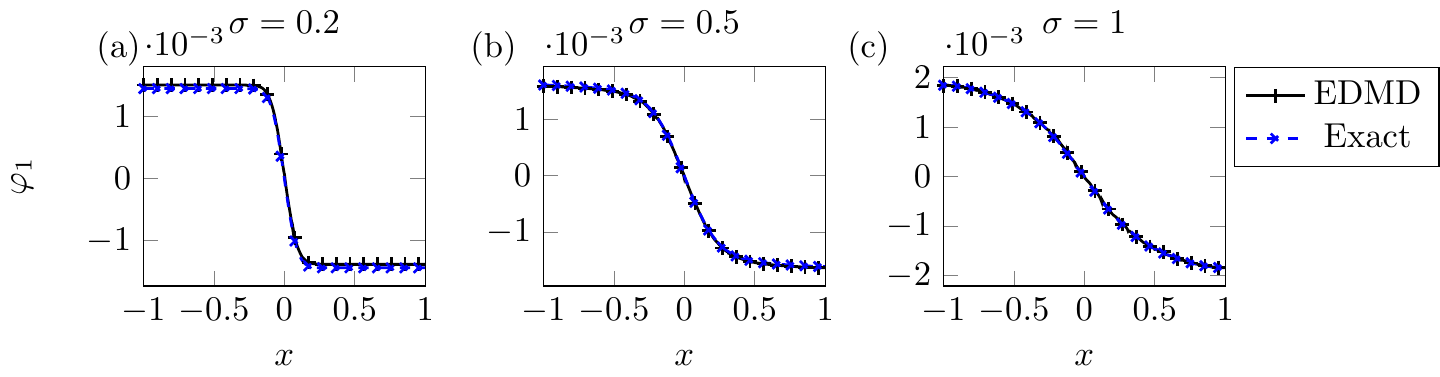}
\caption{A comparison of the leading nontrivial eigenfunction computed with the finite difference method and EDMD for:  (a) $\sigma = 0.2$, (b) $\sigma = 0.5$, and (c) $\sigma = 1.0$. 
As with the eigenvalues, there is excellent agreement between the ``true'' and data driven approximations of this eigenfunction, though there are small quantitative differences in the approximation.
}
\label{fig:double-well-eigenfunctions}
\end{figure}

More interesting is the first {\em nontrivial} eigenfunction  which has the
eigenvalue $\lambda_1$.
The change in $\varphi_1$ as a function of $\sigma$ is shown in Fig.~\ref{fig:double-well-eigenfunctions}.
As with the eigenvalue, there is good agreement between EDMD and the directly computed eigenfunctions at different values of $\sigma$.
For all values of $\sigma$, $\varphi_1$ is an odd function; what changes is how rapidly $\varphi_1$ transitions from its maximum to its minimum. 
When $\sigma$ is small, this transition is rapid, and $\varphi_1$ will approach a step function as $\sigma \to 0$.
When $\sigma$ grows, this eigenfunction is ``smoothed out'' and the transition becomes slower. In the limit as $\sigma\to\infty$, the dynamics of the problem are dominated by the diffusion term, and $\varphi_1$ will be proportional to $\cos(\pi x/L)$ as is implied by the rightmost plot in the figure.
%
%

In many system identification algorithms (e.g., Ref.~\cite{Juang1994}), one often constructs deterministic governing equations from inherently stochastic data (either due to measurement or process noise).
Similarly, methods like DMD have been applied to noisy sets of data to produce an approximations of the Koopman modes and eigenvalues with the assumption that the underlying system is deterministic. 
In this example, this is equivalent to using the output of EDMD with data taken with $0<\sigma \ll 1$ as  an approximation of the Koopman tuples that would be obtained with $\sigma = 0$. 

For certain tuples, this is a reasonable approach. 
Taking $\sigma\to 0$,  $\lambda_3$ and $\lambda_4$  and $\varphi_3$ and $\varphi_4$ are good approximations of their deterministic counterparts. 
In particular $\varphi_3$ and $\varphi_4$ are one-to-one with their associated basin of attraction and appear to possess a zero at the stable fixed point.
However, these approximate eigenfunctions lack some important features such as a singularity at $x =0$ that occurs due to the unstable fixed point there. 
Therefore, both eigenfunctions are good approximations of their $\sigma = 0$ counterparts, but cannot be ``trusted'' in the vicinity of an unstable fixed point.

For other tuples, even a small amounts of noise can be important. 
Consider the ``slowest'' non-zero eigenvalue, $\lambda_2$, which appears to approach $-4$ as $\sigma \to 0$, but is {\em not} obtained by the EDMD method when $\sigma = 0$. 
Formally, the existence of an eigenvalue of $-4$ is not surprising.
The fixed point at $x=0$ is unstable with $\lambda = 4$, and in continuous time, if ($\lambda_n$, $\varphi_n$) is an eigenvalue/eigenfunction pair then ($k\lambda_n$, $\varphi^k_n$) is, at least  formally, an eigenvalue/eigenfunction pair for any scalar $k$.  
Using an argument similar to Ref.~\cite{Matkowsky1981}, it can be shown that $\varphi_2(x) = C_0\exp(-4x^2/\sigma^2) + \mathcal{O}(\sigma^2)$ as $\sigma\to 0$ where $C_0$ is chosen to normalize $\varphi_2$. 
However, this approaches a delta function as $\sigma\to 0$, and therefore leaves the subspace of observables spanned by our dictionary.
When this occurs, this tuple appears to ``vanish,'' which is why it does not appear in the $\sigma = 0$ limit.
As a result, when applying methods like EDMD or DMD to noisy data, the spectrum of the finite dimensional approximation is not necessarily a good approximation of the spectrum that would be obtained with noise--free data. 
Some of the tuples, such as those containing $\varphi_1$, $\varphi_3$, and $\varphi_4$, have eigenvalues that closely approximate the ones found in the deterministic problem. 
However, others such  as the tuple containing $\varphi_2$, do not.
Furthermore, the only method to determine that $\lambda_2$ can be neglected is by directly examining the eigenfunction.
As a result, when we  apply methods like DMD/EDMD to noisy data with the purpose of using the spectrum to determine the time scales and behaviors of the underlying system, we must keep in mind that not all of the eigenvalues obtained with noisy data will be present if ``clean'' data is used instead.
%

%
%

\subsubsection{Rate of Convergence}

Among other things, the performance of the EDMD method is dependent upon the number of snapshots provided to it, the distribution of the data,  the underlying dynamical system, and the dictionary. 
In this section, we examine the convergence of EDMD to a Galerkin method as the number of snapshots increases in order to provide some intuition about the ``usefulness'' of the eigenfunctions obtained without an exhaustive amount of data.
To do so, we generated a larger set of data consisting of $10^7$ initial conditions chosen from a spatially uniform distribution for the case with $\sigma = 1$.
Each initial condition was propagated using the Euler-Maruyama method described in the previous section.
Then we applied EDMD using the same dictionary  to subsets of the data, computed the leading nontrivial eigenvalue and eigenfunction, and compared the results to the ``true'' leading eigenfunction and eigenvalue computed using a finite difference approximation of the stochastic Koopman operator.

\begin{figure}
\includegraphics[width=0.8\textwidth]{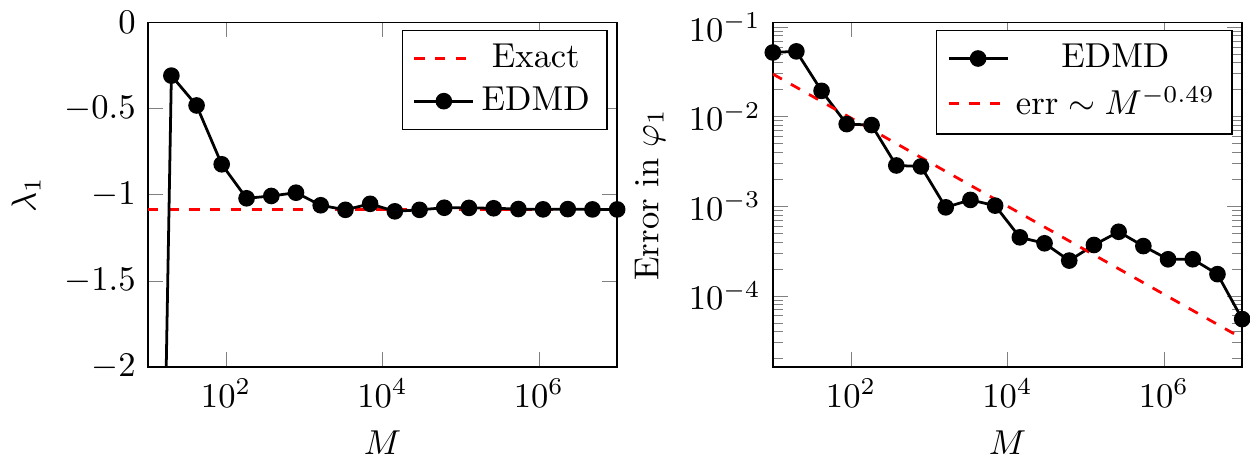}
\caption{(left) Plot of the first nontrivial eigenvalue as a function of the number of data points with $\sigma = 1$. 
The red dashed line denotes the ``exact'' value computed using direct numerical methods. 
Although the EDMD approximation is  poor with $M<100$, it quickly becomes more accurate as $M$ is increased.  
(right) Plot of the error, i.e., $\|\varphi_{1,\text{EDMD}} - \varphi_{1,\text{True}}\|$, as a function of $M$.
Because the scalar products in \eqref{eq:edmd} are evaluated using Monte-Carlo integration, the method converges like $\mathcal{O}(M^{-1/2})$ as shown by the fit on the right.
}
\label{fig:edmd-convergence}
\end{figure}

Figure~\ref{fig:edmd-convergence} shows the convergence of the leading nontrivial eigenvalue and eigenfunction as a function of the number of snapshots, $M$. 
In the rightmost plot, we define the error as $\|\varphi_{1,\text{EDMD}} - \varphi_{1,\text{True}}\|$ after both eigenfunctions have been normalized so that $\|\varphi_{1,\text{EDMD}}\|_2 = \|\varphi_{1,\text{True}}\|_2$.
As  expected, EDMD is inaccurate when $M$ is small (here, $M<100$);  there is not enough data to accurately approximate the scalar products. 
For $M > 10^3$, the eigenfunction produced by EDMD have the right shape, and the eigenvalue is approaching its true value. 
For $M > 10^4$, there is no ``visible'' difference in the leading eigenvalue, and the error in the leading eigenfunction is less than $10^{-3}$.

To quantify the rate of convergence, we fit a line to the plot of error versus $M$ in the right panel of Fig.~\ref{fig:edmd-convergence}. 
As expected, EDMD converges like $M^{-0.49}$, which is very close to the predicted value of $\mathcal{O}(M^{-0.5})$ associated with Monte--Carlo integration. 
Because this problem is stochastic, we cannot increase the rate of convergence by uniform sampling (the integral over the probability space associated with the stochastic dynamics will still converge like $\mathcal{O}(M^{-1/2})$), even though that is a simple method for enhancing the rate of convergence for deterministic problems.

\begin{figure}
\centering
\includegraphics[width=\textwidth]{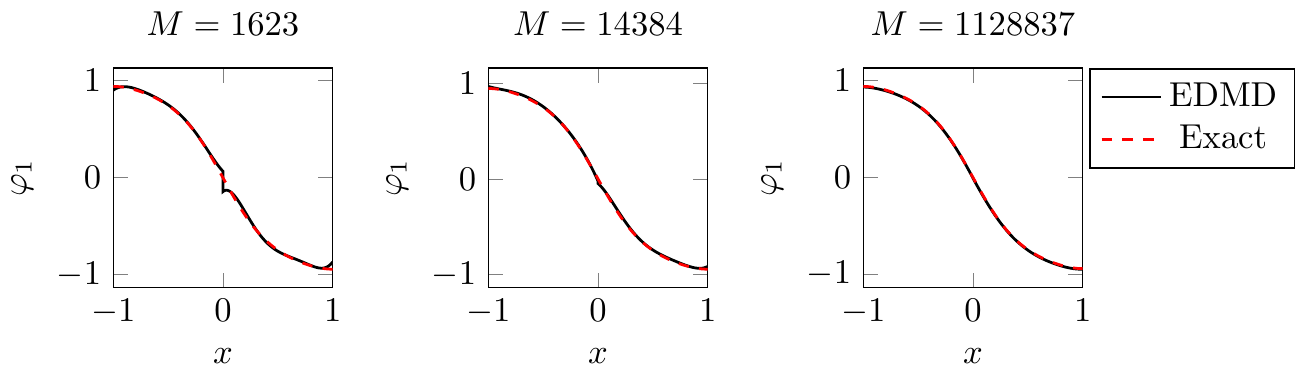}
\caption{Comparison of the EDMD approximation of the leading nontrivial Koopman eigenfunction with $M=1623$, 14384, and 1128837 with the ``exact'' solution obtained using direct numerical methods. 
The EDMD solution is  similar to the exact solution even when $M\sim 10^3$, but we would not characterize it as quantitatively accurate  until $M > 10^4$. 
Beyond $M>10^5$, EDMD and the exact solution are visually identical. 
}
\label{fig:edmd-convergence-eyeball}
\end{figure}

To provide some intuition about what the eigenfunctions look like for a fixed value of $M$, Fig.~\ref{fig:edmd-convergence-eyeball} plots the leading nontrivial eigenfunction for $M=1623$, 14384, and 1128837.
EDMD is qualitatively accurate even at the smallest values of $M$, but there are clearly some numerical issues at the edges of the domain and near $x=0$ where the discontinuities in the numerically computed eigenfunctions can occur with our choice of dictionary.
To obtain a more quantitatively accurate solution, additional data points are required. 
When $M=14384$, the numerical issues at the boundaries and the discontinuity at $x=0$ have diminished.
As shown in the plot with $M = 1128837$, this process continues until the EDMD eigenfunction is visually identical to the true eigenfunction.

\subsection{Parameterizing Nonlinear Manifolds and Reducing Stochastic Dynamics}
\label{subsec:dmaps}

In this section, we will briefly demonstrate how the EDMD method can be used to parameterize nonlinear manifolds and reduce stochastic differential equations defined on those manifolds. 
Everything done here could also be done for a deterministic system; we chose to use an SDE rather than an ODE only to highlight the similarities between EDMD\ and methods like diffusion maps, and not because of any restriction on the Koopman approach.
We proceed in two steps: first, we will show that data from an SDE defined on the Swiss Roll, which is a nonlinear manifold often used as a test of nonlinear manifold learning techniques~\cite{Lee2007, Coifman2006,Nadler2005,Nadler2006}, in conjunction with the EDMD procedure can generate a data driven parameterization of that manifold. 
For this first example, isotropic diffusion is used, so there is no ``fast'' or ``slow'' solution component that can be meaningfully neglected. 
Instead, we will show that the leading eigenfunctions are one-to-one with the ``length'' and ``width'' of the Swiss Roll.
Then we alter the SDE and introduce a ``fast'' component by making the diffusion term anisotropic. 
In this case, EDMD will ``pick'' the slower components before the faster ones.
Both of these tasks can be accomplished using other methods such as Diffusion Maps (DMAPs) and its variants~\cite{Coifman2006,Nadler2005,Nadler2006,Dsilva2013}, and it is only recently that the Koopman operator has also been applied to such problems~\cite{Froyland2013}.

\subsubsection{Parameterizing a Nonlinear Manifold with a Diffusion Process}

For this example, the data are generated by a diffusion process on a rectangular domain,
\begin{equation}
d\vec{s}=2d\vec{W}_{t},
\label{eq:underlying-sde}
\end{equation}
where $\vec{s}=(s_1,s_2)$ is the state and $\vec{W}_{t}$ is a two-dimensional Wiener process with $s_1\in[0,3\pi]$ and $s_2\in[0,2\pi]$.
No flux boundary conditions are imposed at the domain edges.
If one had access to the true variables, the SKO could be written as  
\begin{equation}
\tilde\koopman\phi = 2\partial_{s_1}^2\phi + 2\partial_{s_2}^2\phi,
\label{eq:koopman-dmap}
\end{equation}
also with no-flux boundary conditions; in this particular problem, the SKO is self-adjoint and therefore equivalent to the Perron-Frobenius operator. 
The eigenfunctions should be $\varphi_{ij} = \cos\left(\frac{i}{3}s_1\right) \cos\left(\frac{j}{2}s_2\right)$ with the eigenvalues $\lambda_{ij} = -2\left(\frac{i^2}{9}+\frac{j^2}{4}\right)$.
Note that the leading eigenfunctions, $\varphi_{1,0}$ and $\varphi_{0,1}$ are $\cos\left(\frac{1}{3}s_1\right)$ and $\cos\left(\frac{1}{2}s_2\right)$, which are one-to-one with $s_1$ and $s_2$ on $[0, 3\pi]$ and $[0, 2\pi]$ respectively, and could be used to parameterize state space if $s_1$ and $s_2$ were not known.

In this example, these true data (i.e., the state expressed in terms of $s_1$ and $s_2$) on the rectangle are mapped onto a ``Swiss Roll'' via the transformation
\begin{equation}
\vec{g}(\vec{s)}=\begin{bmatrix}(s_1+0.1)\cos(s_1)\\
s_2\\
(s_1+0.1)\sin(s_1)
\end{bmatrix},
\label{eq:swiss-mapping}
\end{equation}
which, among other things, has introduced a new spatial dimension. 
In all that follows, {\em the EDMD approach is applied to the 3-dimensional, transformed variables and not the 2-dimensional, true variables}.
Our objective here is to determine a 2-parameter description of what initially appears to be 3-dimensional data, directly from the data.

\begin{figure}
\centering
\includegraphics[width=\textwidth]{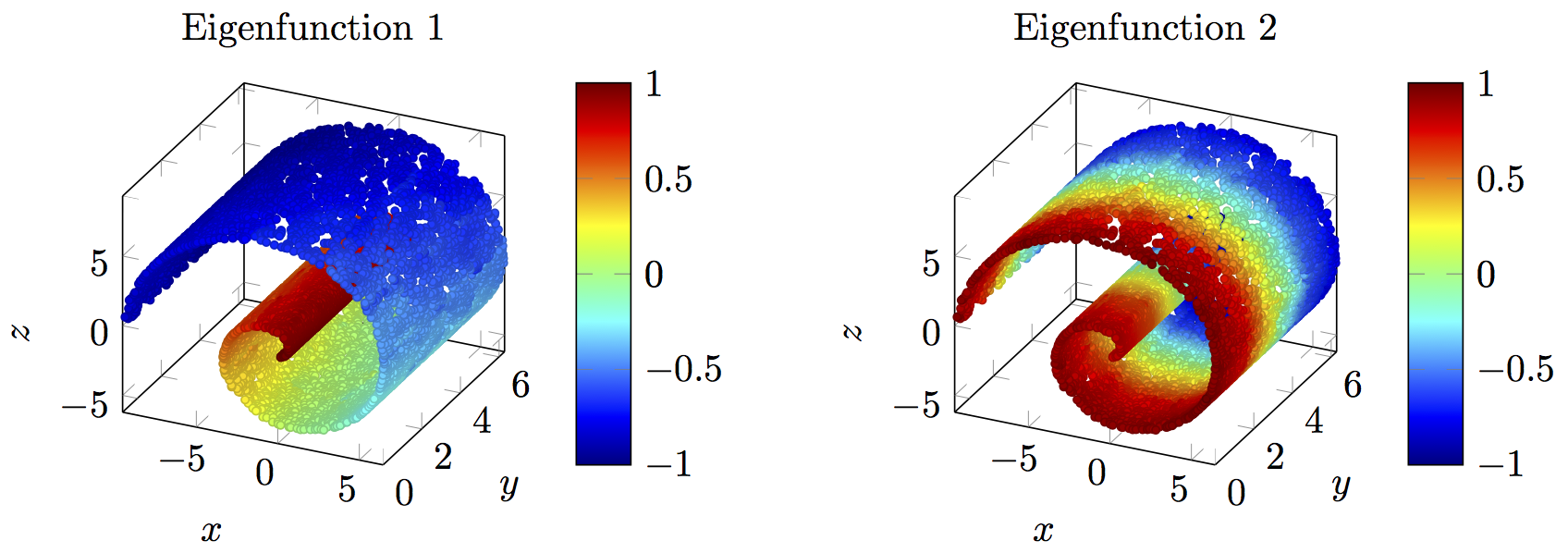}
\caption{The first two nontrivial Koopman eigenfunctions for the diffusion process on the ``Swiss Roll'' using $3\times 10^4$ data points.  
The first eigenfunction is one-to-one with $s_1$ (the ``length'' of the Swiss Roll), and the second eigenfunction is one-to-one with $s_2$.
As a result, they could act as a data driven parameterization of this nonlinear manifold. 
The eigenvalues associated with these eigenfunctions are -.234 and -0.491 compared to the theoretical values of $-\frac{2}{9}$ and $-0.5$ respectively. 
}
\label{fig:diffusion-map}
\end{figure}

The data given to EDMD were generated by $10^4$ initial conditions uniformly distributed in $s_1$ and $s_2$ that were evolved for a total time of $\Delta t = 0.1$ using the Euler-Maruyama method with 100 timesteps.
Then both the initial and terminal states of the system were mapped into 3 dimensions using \eqref{eq:swiss-mapping}. 
Next, a dictionary must be defined. 
However, $\manifold$ is unknown (indeed, parameterizing $\manifold$ is the entire point), so $\Omega$ is embedded in $\mathbb{R}^3$ such that $x\in[-3\pi-.1, 3\pi+.1]$, $y\in[0, 2\pi]$ and $z\in[-3\pi-.1, 3\pi+.1]$.
In this larger domain, the spectral element basis consisting of 4096
rectangular subdomains (16 each in $x$, $y$, and $z$) with up to linear polynomials in each subdomain
is employed. 
Because $\manifold\subset \Omega$, extraneous and redundant trial functions are expected, and  $\mat{G}$ is often ill conditioned.

Figure~\ref{fig:diffusion-map} shows the transformed data colored by the first and second nontrivial eigenfunctions of the Koopman operator.
Unlike many of the previous examples, there are clear differences in the analytical and computed eigenvalues and eigenfunctions.
However, the first eigenfunction is one-to-one with the ``arclength'' along the Swiss roll (i.e., the $s_1$ direction), and the second eigenfunction is one-to-one with the ``width'' (i.e., the $s_2$ direction).
Furthermore, the first two eigenvalues obtained with EDMD, -.234 and -0.491, arrange the eigenfunctions in the correct order and compare favorably with the true eigenvalues of $-2/9$ and $-1/2$. 
Therefore, while our computation of the Koopman eigenfunctions may not be highly accurate, they are  accurate enough to parameterize the nonlinear manifold, which was the goal. 

The procedure for incorporating new data points is simple; the embedding for any $\tilde{\vec x}\in\manifold$ can be obtained simply by evaluating the relevant eigenfunctions at $\tilde{\vec x}$.  
It should be stressed that although the $\varphi$ are defined on $\Omega$, their value is only meaningful on (or very near) $\manifold$ because that is where the dynamical system is defined. 
Therefore, these new points must be elements of $\manifold$ if the resulting embedding is to have any meaning.

\subsubsection{Reducing Multiscale Dynamics}

In the previous example, the noise was isotropic, so the dynamics were equally ``fast'' in both directions. 
As a result, the $s_1$ component was prioritized because the underlying rectangular domain is larger in $s_1$ than it is in $s_2$. 
In this example, we introduce anisotropic diffusion, and therefore create  ``fast'' and ``slow'' directions on the nonlinear manifold. 
The purpose of this example is to show that EDMD will ``reorder'' its ranking of the eigenfunctions, and recover the slower component before the faster one if the level of anisotropy is large enough.

Our particular example is
\begin{subequations}
\begin{align}
ds_1 &= \frac{2}{\epsilon}dW_1, \\
ds_2 &= 2dW_2,
\end{align}
\label{eq:multiscale-sde}
\end{subequations}
with $\epsilon=0.1$, which is again transformed onto the same Swiss Roll.
Although the domain in $s_1$ is larger than it is in the $s_2$ component, the dynamics of $s_1$ are now significantly faster than those of $s_2$ due to a much larger amount of ``noise''. 
Because the two random processes in \eqref{eq:multiscale-sde} are independent, the eigenfunctions themselves should not change. 
However, the {\em eigenvalues} associated with each of the eigenfunctions should.

\begin{figure}
\centering
\includegraphics[width=\textwidth]{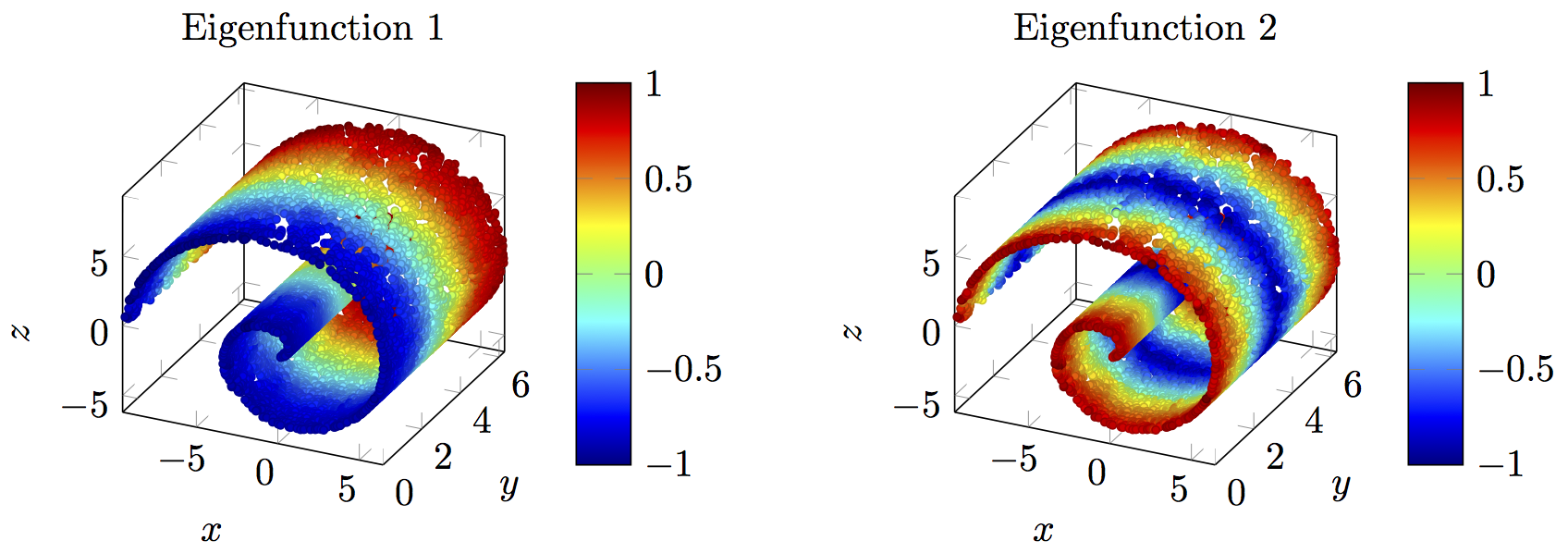}

\caption{The first two nontrivial Koopman eigenfunctions for the multiscale diffusion process in \eqref{eq:multiscale-sde} on the ``Swiss Roll'' computed using $3\times 10^4$ data points.  
The eigenvalues associated with these eigenfunctions are -0.504 (analytically it has the value -1/2) and -2.1 (analytically, -2).
In contrast with Fig.~\ref{fig:diffusion-map}, these eigenfunctions are a parameterization of the ``width'' of the Swiss roll ($s_2$) and a higher charmonic. 
As a result, the eigenfunctions computed by EDMD take into account both the geometry of the underlying manifold and the dynamics defined on it.
}
\label{fig:swiss-noisy}
\end{figure}

Figure~\ref{fig:swiss-noisy} shows how EDMD captures this difference in the underlying diffusion process.
Before, the first two nontrivial eigenfunctions were one-to-one with $s_1$ and $s_2$ respectively; now, the first is one-to-one with $s_2$, and the second is a higher harmonic  (but still only a function of $s_2$). 
The eigenfunction that is one-to-one with $s_1$ still exists, but it is no longer associated with a leading eigenvalue.  
Analytically, its eigenvalue is $-20/9$ though EDMD computes a value of $-2.3$.
Therefore, it is no longer a leading eigenfunction, but still is approximated by the EDMD procedure. 

This section explored the application of the EDMD method to data taken from a Markov process. 
Algorithmically, the method remains the same regardless of how the data were generated, but as demonstrated here, EDMD computes an approximation of the tuples associated with SKO rather than the Koopman operator.
To demonstrate the effectiveness of EDMD, we applied the method to a simple SDE with a double--well potential, and an SDE defined on a Swiss Roll, which is a nonlinear manifold often used as a benchmark for manifold learning techniques.
One advantage of the Koopman approach for applications such as manifold learning or model reduction is that the Koopman tuples take into account both the geometry of the manifold, through the eigenfunction and mode, and the dynamics, through the eigenvalue.   
As a result, the approach taken here is aware of both geometry and dynamics, and does not focus solely on one or the other.

\section{Conclusions}
\label{sec:conclusions}

In this manuscript, we presented a data driven method that computes approximations of the Koopman eigenvalues, eigenfunctions, and modes (what we call Koopman tuples) directly from a set of snapshot pairs.
We refer to this method as Extended Dynamic Mode Decomposition (EDMD). 
The finite dimensional approximation generated by EDMD is the solution to a least squares problem, and converges to a Galerkin method with a large amount of data. 
While the usefulness of the Galerkin method depends on the sampling density and dictionary selected, several ``common sense'' choices of both appear to produce useful results.

We demonstrate the effectiveness of the method with four examples: two examples dealt with deterministic data, and two with stochastic data. 
First, we applied EDMD to a linear system where the Koopman eigenfunctions are known analytically. 
Direct comparison of the EDMD eigenfunctions and the analytic values demonstrated that EDMD can be highly accurate with the proper choice of data and  dictionary. 
Next, we applied EDMD to the unforced Duffing oscillator where the Koopman eigenfunctions are not known explicitly.
Although more data will increase the accuracy of the resulting eigenfunctions, they appeared to be accurate enough to effectively partition the domain of interest and parameterize the resulting partitions.

The final two examples used data generated by Markov processes.
First, we applied EDMD to data taken from an SDE with a double well potential, and demonstrated the accuracy of the method by comparing those results with a direct numerical approximation of the stochastic Koopman operator over a range of diffusion parameters.
Next, we applied EDMD to data from a diffusion process on a ``Swiss Roll,'' which is a nonlinear manifold commonly used as an example for nonlinear dimensionality reduction. 
Similar to those methods (see e.g., Ref.~\cite{Coifman2006, Lee2007}), EDMD generated an effective parameterization of the manifold using the leading eigenfunctions. 
By making the diffusion anisotropic, we then demonstrated that EDMD extracts a parameterization that is dynamically, rather than  only  geometrically, meaningful.
Due to the simplicity of this problem, the eigenfunctions remain unchanged despite the anisotropy; the difference appears in the temporal evolution of the eigenfunctions, which is dictated by the corresponding set of eigenvalues.  
As a result, the purpose of that example was to show that EDMD ``ordered'' the eigenvalues of each tuple appropriately.

The Koopman operator governs the evolution of observables defined on the state space of a dynamical system.
By judiciously selecting how we observe our system, we can generate {\em linear models} that are valid on  all of (or, at least, a larger  subset of) state space rather than just some small neighborhood of a fixed point; this could allow algorithms designed for linear systems to be applied even in nonlinear settings.
However, the tuples of eigenvalues, eigenfunctions, and modes required to do so are decidedly nontrivial to compute. 
Data driven methods, such as EDMD, have the potential to allow accurate approximations of these quantities to be computed without knowledge of the underlying dynamics or geometry.
As a result, they could be a practical method for enabling Koopman-based analysis and model reduction in large nonlinear systems.

\section*{Acknowledgments}

The authors would like to thank Igor Mezi{\'c}, Jonathan Tu, Maziar Hemati, and Scott Dawson for interesting and useful discussions on Dynamic Mode Decomposition and the Koopman operator. 
M.O.W. gratefully acknowledge support from NSF DMS-1204783.  
I.G.K acknowledges support from AFOSR FA95550-12-1-0332 and NSF CMMI-1310173.
C.W.R acknowledges support from AFOSR FA9550-12-1-0075.

\appendix
\section{EDMD with Redundant Dictionaries}
\label{app:redundant}

In this appendix, we present a simple example of applying EDMD to a problem where the elements of $\dictionary$ contain redundancies (i.e., the elements of $\dictionary$ are not a basis for  $\observables_{\dictionary} \subset \observables$).   
Given full knowledge of the underlying dynamical system, one would always choose the elements of $\dictionary$ to be a basis for $\observables_{\dictionary}$, but due to our ignorance of $\manifold$, a redundant set of functions may be chosen.
Our objective here is to demonstrate that accurate results can still be obtained even if such a choice is made.
To separate quadrature errors from errors resulting from our choice of $\dictionary$, we assume that $M$ is large enough that the EDMD method has already converged to a Galerkin method in that the residual is orthogonal to the space spanned by $\dictionary$.

For the purposes of demonstration, we replace $\koopman$ with $\mathcal{L}=\partial_s^2$, the Laplace--Beltrami operator defined on the manifold, $\manifold$, where $(x,y) = (s, s)$ for $s\in[0, 2\pi)$ with periodic boundary conditions, which would correspond to, say, the EDMD procedure applied to a diffusion process on a periodic domain. 
A useful basis for this problem would be $\tilde\psi_k(x,y) = \exp(\imath ks) = \exp(\imath k(x+y))$, but without prior knowledge of $\manifold$, it is difficult to determine this choice should be made.
Because the problem appears two dimensional, one may choose a dictionary whose elements have the form  $\psi_{m, n}(x, y) = \exp(\imath m x + \imath ny)$, which contains the $\tilde\phi$ but is not linearly independent on $\manifold$. 
The indexes we use for the trial functions are $\psi_k(x,y) = \psi_{m, n}(x,y)$ with  $m = (k \mod K) - K/2$ and $n = \left\lfloor\frac{k}{K}\right\rfloor - K/2$ with $k = 0, 1, \ldots, K^2$.
Here $K\in\mathbb{N}$ is the total number of basis functions in a single spatial dimension.

\begin{figure}
\includegraphics[width=\textwidth]{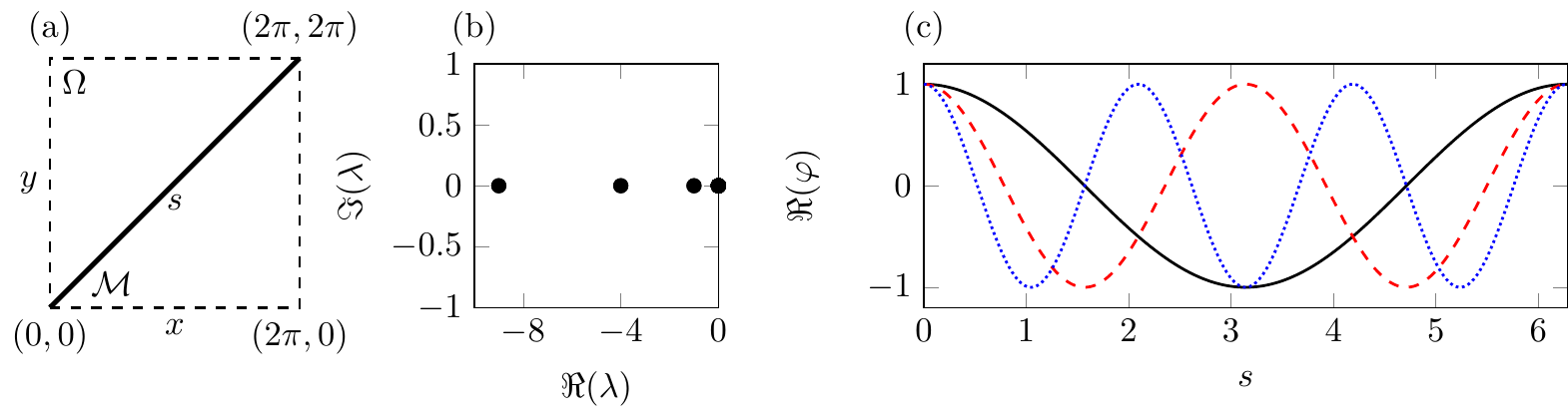}
\caption{
(a) A sketch of the manifold $s\mapsto (s, s)$ where our dynamical system is defined, and the larger domain, $\Omega$, on which the elements of $\dictionary$ are defined. 
(b)  A plot of the leading 56 eigenvalues of $\partial_s^2$ computed using EDMD; the redundant trial functions have increased the dimension of the nullspace from 1 to 50, but accurately capture the pairs of eigenvalues at $-k^2$ for $k=0,1,2,\ldots,8$.
(c) A plot of the real part of the first three nontrivial eigenfunctions shown in black, red, and blue respectively; as expected, they are equivalent to $\cos(ks)$.
The imaginary component of the eigenfunctions, which is not shown, captures the $\sin(ks)$ terms.  
}
\label{fig:edmd-redundant}
\end{figure}

Following \eqref{eq:edmd}, the $i,j$-th element of $\mat{G}$ is 
\begin{equation}
\begin{aligned}
\mat{G}_{i,j} &= \int_{\manifold} \psi_i(\vec x)^* \psi_j(\vec x) d\vec x = 
\int_{0}^{2\pi}  e^{\imath ((m_j - m_i)s + (n_j - n_i)s)} \;ds \\
&= 
\begin{cases}
2\pi & m_j + n_j  - m_i - n_i = 0, \\
0    & \text{otherwise}.
\end{cases}
\end{aligned}
\end{equation}
Similarly, 
\begin{equation}
\begin{aligned}
{\mat{A}}_{i,j} &= \int_{\manifold} \psi_i(\vec x)^* \partial_s^2\psi_j(\vec x) d\vec x = 
\int_{0}^{2\pi}  -(m_j + n_j)^2e^{\imath ((m_j - m_i)s + (n_j - n_i)s)} \;ds \\
&= 
\begin{cases}
-2\pi(m_j + n_j)^2 & m_j + n_j  - m_i - n_i = 0, \\
0    & \text{otherwise}.
\end{cases}
\end{aligned}
\end{equation}
The diagonal structure we would normally have has been replaced with a more complex sparsity pattern, and it has a large nullspace (when $K=8$, the nullspace is 50 dimensional). 
To reiterate, there are no advantages to this choice; the redundancies in $\dictionary$ appear due to ignorance about the nature of $\manifold$, which is the expected situation. 
Because $\mat{G}$ is singular, the use of the pseudoinverse in \eqref{eq:edmd} is {\em critical} to obtain a unique solution.

However, once this is done, there is excellent agreement between the leading eigenfunctions and eigenvalues of $\mathcal{L}$ and those computed using EDMD; this is shown in Fig.~\ref{fig:edmd-redundant}.
The nonzero eigenvalues are quantitatively correct; in particular, pairs of eigenvalues of the form $\lambda=-k^2$ are obtained up until $k=8$ using $K = 8$.
Although the maximum (absolute) value of $m$ or $n$ is only 4, it can be seen from the form of the trial functions that the superposition of these functions on $\manifold$ mimics $k = 8$ modes.
The associated eigenfunctions are shown in Fig.~\ref{fig:edmd-redundant}c; again, there is excellent agreement between the analytic solution (i.e., $\exp(-\imath ks)$) and the EDMD computed solution.

The resulting eigenfunctions can also be evaluated for $(x,y)\not\in\manifold$, but the functions have no dynamical meaning there.
Indeed, their value is determined entirely by the regularization used, and has no relationship to the underlying dynamical system, which is defined solely on $\manifold$.
This should be contrasted to related works such as Ref.~\cite{Froyland2013} where the dynamical system is truly defined on $\Omega$, and $\manifold$ is simply the slow manifold where the eigenfunctions evaluated at $(x,y)\not\in\manifold$ are meaningful as they contain information about the fast dynamics of the system.

Overall, the performance of the EDMD procedure is dependent upon the subspace, $\observables_{\dictionary}$ and not the precise choice of $\dictionary$.
There are numerical advantages to choosing $\dictionary$ to be a basis for $\observables_{\dictionary}$, but in many circumstances this cannot be done without prior knowledge of $\manifold$.
As a result, there are likely benefits to combining EDMD with manifold learning techniques (see, e.g., Refs.~\cite{Coifman2006, Lee2007,Erban2007,Sirisup2005}).  
These methods can numerically approximate $\manifold$, which could allow a more effective choice of the elements of $\dictionary$ and their associated numerical benefits.
As shown here, these methods are not essential to the algorithm, but $\manifold$ must be identified through some means if EDMD is to be used for more than just data analysis.

\bibliographystyle{ieeetr}
\bibliography{koopman}

\end{document}